\documentclass[11pt,a4paper]{article}
\usepackage[utf8]{inputenc}
\usepackage[english]{babel}
\usepackage[T1]{fontenc}
\usepackage{amscd,amssymb,amsmath,amsthm,amstext,amsfonts}
\usepackage[amssymb]{SIunits}      % measure symbols 
\usepackage{stmaryrd}
\usepackage[pdftex]{graphicx}
\usepackage{fancybox}
\usepackage[square,numbers]{natbib}
\usepackage{booktabs}
\usepackage{hyperref}
\usepackage{algorithm}
\usepackage{algpseudocode}
\usepackage{enumerate}  
\usepackage[capitalise,nameinlink]{cleveref}
\crefformat{equation}{\textup{#2(#1)#3}}
\crefname{figure}{Figure}{Figures}
\usepackage{subfigure}
\usepackage{comment}
\usepackage{mathtools}
\usepackage{tikz,pgfplots}
\pgfplotsset{compat=newest}
\usepackage{pgfplotstable}
\pgfplotsset{compat=newest}
\usepackage{multirow}
\usepackage{authblk} 
\usepackage[normalem]{ulem}
\usepackage{amsopn}

\usepackage{listings}
\usepackage{color}
\usepackage{todonotes}
\usepackage{bm}

\theoremstyle{definition}

\theoremstyle{remark}

\numberwithin{theorem}{section}
\numberwithin{equation}{section}
\numberwithin{table}{section}
\numberwithin{figure}{section}

\usepackage{bm}
\newcommand{\ee}{\bm{e}}

\newcommand{\nn}{\bm{n}}
\newcommand{\uu}{\bm{u}}
\newcommand{\ii}{\bm{i}}
\newcommand{\jj}{\bm{j}}
\newcommand{\pp}{\bm{p}}
\newcommand{\qq}{\bm{q}}
\renewcommand{\tt}{\bm{t}}
\newcommand{\qoi}{f}
\newcommand{\Ev}{\mathbb{E}}
\newcommand{\Nset}{\mathbb{N}}
\newcommand{\Rset}{\mathbb{R}}
\newcommand{\Var}{\mathbb{V}}
\newcommand{\GGG}{\mathcal{G}}
\newcommand{\II}{\mathcal{I}}
\newcommand{\SG}{\mathcal{S}}
\newcommand{\degr}{r}
\newcommand{\degq}{q}
\newcommand{\ncp}{N_{coll}}
\newcommand{\mcp}{M_{coll}}

\allowdisplaybreaks

\author[1]{Giuseppe Balduzzi}
\affil[1]{DICAr - Department of Civil Engineering and Architecture, University of Pavia, Italy}

\author[2]{Francesca Bonizzoni}
\affil[2]{MOX - Dipartimento di Matematica, Politecnico di Milano, Italy}

\author[3]{Lorenzo Tamellini}
\affil[3]{ Istituto di Matematica Applicata e Tecnologie Informatiche ``E. Magenes'' - Consiglio Nazionale delle Ricerche (IMATI-CNR), Pavia, Italy}

\title{Uncertainty quantification in timber-like beams using sparse grids: theory and examples with off-the-shelf software utilization}

\date{}

\begin{document}

\definecolor{mygreen}{RGB}{28,172,0} % color values Red, Green, Blue
\definecolor{mylilas}{RGB}{170,55,241}
\definecolor{mylightgray}{RGB}{220,220,220} % a shade lighter than standard lightgray
% base font for code
\def\lstbasicfont{\bfseries\fontfamily{pcr}\selectfont\footnotesize}

\lstset{language=Matlab,%
    %basicstyle=\color{red},
    basicstyle={\lstbasicfont},
    backgroundcolor = \color{mylightgray}, %background color of code
    captionpos = b, % sets the captions BELOW the code
    frame=single,rulecolor=\color{gray}, % frames the code
    breaklines=true,%
    morekeywords={matlab2tikz},
    keywordstyle=\color{blue},%
    morekeywords=[2]{1}, keywordstyle=[2]{\color{black}},
    identifierstyle=\color{black},%
    stringstyle=\color{mylilas},
    commentstyle=\color{mygreen},%
    showstringspaces=false,%without this there will be a symbol in the places where there is a space
    numbers=left,%
    numberstyle={\tiny \color{black}},% size of the numbers
    numbersep=9pt, % this defines how far the numbers are from the text
    emph=[1]{for,end,break},emphstyle=[1]\color{red}, %some words to emphasise
    %emph=[2]{word1,word2}, emphstyle=[2]{style},    
}

\maketitle

\begin{abstract}
When dealing with timber structures, the characteristic strength and stiffness of the material are made highly variable and uncertain by the unavoidable, yet hardly predictable, presence of knots and other defects.
In this work we apply the sparse grids stochastic collocation method to perform uncertainty quantification for structural engineering in the scenario described above.
Sparse grids have been developed by the mathematical community in the last decades and their theoretical background has been rigorously and extensively studied.
The document proposes a brief practice-oriented introduction with minimal theoretical background, provides detailed instructions for the use of the already implemented Sparse Grid Matlab kit (freely available on-line) and discusses two numerical examples inspired from timber engineering problems that highlight how sparse grids exhibit superior performances compared to the plain Monte Carlo method.
The Sparse Grid Matlab kit requires only a few lines of code to be interfaced with any numerical solver for mechanical problems (in this work we used an isogeometric collocation method) and provides outputs that can be easily interpreted and used in the engineering practice.
\end{abstract}

\vspace{1cm}

% % % % % % % % % % % % % % % % % % % % 
\section{Introduction}
% % % % % % % % % % % % % % % % % % % % 

Timber is one of the oldest building material. 
Used since the prehistory, wood has been employed in all ages and by all civilizations, often with peculiar technologies \citep{ss_08}.
Between the Nineteenth and Twentieth centuries other materials (like cast iron, steel, alluminium, and concrete) became largely available, deeply impacting word economic development and sustaining human expansion \citep{nptm_05}.
In recent years, climate change emerged as a new, urgent problem and construction and related industries (in particular, concrete and steel ones) are  the ones with greatest environmental impact \citep{xwz_18, tj_06}.
In this context, timber and wood-based structural elements are experiencing a new springtime.
Indeed, wood is a renewable resource (if forests and production processes are properly managed) \citep{gcc_21},
it presents extremely high strength vs weight ratio, low heat conductivity, good durability, and fire resistance when suitable expedients are implemented \citep{bs_17}.
Furthermore, photosynthesis traps a significant amount of carbon dioxide ($\approx 50 \%$ of dry wood is constituted by carbon) that remains within wood for its whole life \citep{gs_11}.
As a consequence, it represents the best candidate for replacing materials with more significant environmental impact.

Unfortunately, natural growth and sawing of logs lead to a huge variability 
and uncertainty on the mechanical properties of wood: \citep{bs_17} specifies that the strength of wood specimens can change by an order of magnitude, even within the same wood species.
In particular, knots - resulting from the insertion of branches in the stem - often coincide with the point where cracks start, representing therefore the weak point of structural elements \citep{vlef_22}.
Such a situation does not allow an economically convenient exploitation of the material and several strategies have been developed for limiting the negative influence of defects on the performance of the structural element.
The oldest is represented by grading, that consists in different procedures and technologies aiming at sorting sawn timber (and boards used for the manufacturing of glued laminated timber beams and cross laminated timber plates) in classes with assigned characteristic strength \citep[Article B5]{bs_17}.
Nowadays, novel technologies - like laser scanners \citep{p_10} and X-ray computer tomography \citep{lm_ea_12} - allow for the detection of grain direction and wood density, which have been employed for the reconstruction of knot geometry and estimation of wood mechanical properties \citep{klf_16, lkhof_19}.
However, the evaluation of mechanical behavior of wood is characterized by high levels of uncertainty, despite the continuous development of analysis and manufacturing technologies:
a quantitative assessment of how the uncertainty on the mechanical behavior of wood translates to uncertainty on the structural behavior of timber construction can be done
by means of Uncertainty Quantification (UQ) techniques.

UQ techniques applied to timber structures have been object of preliminary investigations in the recent engineering literature
by a multitude of approaches ranging from standard sampling methods as in \citep{jlgk_15} to more advanced methods such as the Stochastic Galerkin, the method of moment equations and Gaussian Processes \citep{kfe_15, fke_16, kf_17, cskd_22}. However, UQ approaches as discussed in the engineering literature are usually targeted on specific engineering problems. 
Furthermore, they often propose imprecise treatments of the mathematical aspects of the problem (like, e.g., error estimates and convergence rates), typically preventing the immediate comparison with other available methods and, ultimately, the choice of the most performing one.
On the contrary, mathematicians have been developing the above-mentioned and several other methods as general-purpose tools, providing also detailed theoretical results; in particular, see e.g. \cite{lk_10,eigel:adaptive} for a general introductions to Stochastic Galerkin, \cite{bn_20,bnk_16,bn_14,bbn_13} for the method of moment equations, and \citep{rasmussen:gpbook2006} for Gaussian Processes. Unfortunately, such methods are often complex to be implemented and might not be available as ready-to-use software, discouraging practitioners from their use.

The present contribution deals with the numerical discretization of the equilibrium equations under uncertainty for a timber-like planar body:
more precisely, we assume that its elasticity modulus depends on a set of parameters modeling the random location and shape of knots, under the simplistic assumptions of heterogeneous and isotropic material. As a result, the displacement vector depends on both the space variable and the set of parameters. The equations - along the space variable - are discretized following the isogeometric analysis (IGA) principles, specifically the IGA collocation that combines high-performance with easy implementation, thanks to the possibility of directly using the strong formulation of the problem. The parametric dependence is treated using the stochastic collocation method based on Smolyak sparse grids, an efficient UQ technique proposed and deeply analyzed by the mathematical community during the latest decades \cite{babuska.nobile.eal:stochastic2,xiu.hesthaven:high}
and implemented in several packages like, e.g., the Sparse Grid Matlab kit \cite{piazzola.tamellini:SGK}. The sparse grids methdology is most effective when the problem at hand depends on a moderate number of uncertain parameters (say up to 20/30 parameters, even though applications to problem with hundreds of random variables are available in literature \cite{chen:adaptive.lognormal,ernst.eal:collocation-logn}), and the outputs
of the model depend smoothly on the input parameters. 

The main contributions of the present work are: (i) the superiority of the stochastic collocation method with respect to the plain Monte Carlo method is demonstrated by means of several numerical tests in the continuum mechanics framework; (ii) algorithmic and implementation details are provided to show its ease of use and possible application to any structural engineering problem.

The rest of the paper is organized as follows. In Section~\ref{sec:model_problem} we introduce the problem of interest, namely the elasticity equation for timber-like beams, where the material variability is encoded in a set a parameters;
moreover, Section~\ref{sec:IGA} details the numerical scheme applied to discretize the model problem in the physical variable. Section~\ref{sec:sparse_grids} is dedicated to the UQ methodology that we employ throughout the work. In Section~\ref{sec:numerical_experiments} a forward UQ analysis is performed on two numerical experiments, namely the expectation and the probability density function of selected quantities of interest are computed and the global sensitivity analysis is carried out.
The conclusions are finally driven in Section~\ref{sec:conclusions}.

% % % % % % % % % % % % % % % % % % % % 
\section{Deterministic mechanical problem}
\label{sec:model_problem}
% % % % % % % % % % % % % % % % % % % % 

\subsection{Continuum mechanic PDEs} \label{sec:pdes}

Let $D=[0,L]\times[0,H]\subset\mathbb{R}^2$ denote a two-dimensional timber beam with length $L>0$ and height $H>0$. Let $\pmb{\mathbb{C}}$ denote the fourth order stiffness tensor, which is assumed to depend on the space variable $(x,y)\in D$ as well as on a set of $N$ parameters $\pp=(\pp_1,\ldots,\pp_N)$ randomly varying in the hyperrectangle $\Gamma\coloneqq\Gamma_1\times\cdots\times\Gamma_N\subset\Rset^N$, with $\Gamma_n=[a_n,b_n]\subset\Rset$ for all $n=1,\ldots,N$.
In particular, $\pmb{\mathbb{C}}$ assumes the following form: 
\begin{equation}
    \label{eq:C}
    \pmb{\mathbb C}(x,y,\pp)
    =
    \begin{bmatrix}
	    E(x,y,\pp) & 0 & 0\\
	    0 & E(x,y,\pp) & 0\\
	    0 & 0 & \frac{E(x,y,\pp)}{2}
    \end{bmatrix}, 
\end{equation}
where the (positive) parameter-dependent elasticity modulus $E(x,y,\pp)$ is modeled as
\begin{equation}
    \label{eq:E}
    E(x,y,\pp)=E_0\,\alpha(x,y,\pp).
\end{equation}
Specific information on the value assumed by $E_0\in\Rset_+$ as well as the form of the function $\alpha(x,y,\pp): D \rightarrow \Rset_+$ will be provided in Section~\ref{sec:numerical_experiments}. 

Given a parameter-independent external load $\tt=(t_x,t_y)$, we look for the displacement $\uu=(u_x,u_y)\colon D\times\Gamma\rightarrow\mathbb{R}^2$ such that 
\begin{equation}
    \label{eq:pde}
    \left\{\begin{array}{ll}
         \operatorname{div}(\pmb{\mathbb{C}}(x,y,\pp)\, \colon \nabla^s \uu(x,y,\pp)) = \mathbf{0}, & (x,y)\in D,\\
        (\pmb{\mathbb{C}}(x,y,\pp)\, \colon \nabla^s \uu(x,y,\pp))\cdot\nn = \tt(x,y), & (x,y)\in\Sigma_t,\\
        \uu(x,y,\pp) = \bm{0}, & (x,y)\in\Sigma_s,
    \end{array}\right.
\end{equation}
where $\{\Sigma_t,\Sigma_s\}$ is a partition of $\partial D$ and $\nabla^s$ denotes the symmetric gradient. The differential operators in~\eqref{eq:pde} are intended with respect to the physical variables $x,\, y$. 
Note that, in the present paper, the beam material is assumed heterogeneous (since the stiffness tensor depends on $x,y$) and isotropic. The latter assumption is not fulfilled in the specif case of timber beams. Nonetheless, it simplifies the theoretical and numerical treatment of the addressed problem. The generalization of the presented results to the anisotropic framework is worth investigating, and will be addressed in a future contribution.

% % % % % % % % % % % % % % % % % % % % 
\subsection{IGA discretization in the space variables}
\label{sec:IGA}
% % % % % % % % % % % % % % % % % % % % 

Using the notation on provided in \ref{appendix:basic_b-splines}, we look for approximations to $u_x, u_y$ of the form
\begin{equation*}
\label{eq_approximation}
    \begin{aligned}
    u_x(x,y)&\approx \sum_{i=1}^{\ncp}\sum_{j=1}^{\mcp} 
    (\hat{u_x})_{i,j}
    R_{i,j}^{\degr,\degq} \left( x,y \right)\\
    u_y(x,y)&\approx \sum_{i=1}^{\ncp}\sum_{j=1}^{\mcp} 
    (\hat{u_y})_{i,j}
    R_{i,j}^{\degr,\degq} \left( x,y \right) 
    \end{aligned}
\end{equation*}
where $R_{i,j}^{\degr,\degq} \left( x,y \right)$ are
bi-variate B-splines, and we require them to be strong solutions to Equation \eqref{eq:pde}.

The obtained equations are then collocated at the Greville abscissae $\left( \hat{x}_i , \hat{y}_j \right)$ ($i=1,\dots,\ncp-1$, $j=1,\dots,\mcp-1$), which can be computed as:
\begin{equation}
\begin{aligned}
  \label{eq:greville}
  \hat{x}_i = \frac{x_{i+2} + x_{i+3} + \ldots + x_{i+r}}{r - 1} , \ i=1,\dots,\ncp-1,
\\
  \hat{y}_j = \frac{y_{j+2} + y_{j+3} + \ldots + y_{j+q}}{q - 1} , \ j=1,\dots,\mcp-1.
\end{aligned}
\end{equation}

The resulting algebraic system of equations, consisting of $2\left(\ncp-1\right)\left(\mcp-1\right)$ equations in the $2\ncp\mcp$ unknowns (i.e., $\ncp \times \mcp$ unknowns for both $u_x$ and $u_y$), has to be finally completed by $2\ncp + 2\mcp -2$ suitable boundary conditions to be imposed as additional equations, as specified in \cite{abhrs_12}.

% % % % % % % % % % % % % % % % % % % % 
\section{Sparse grids and Uncertainty Quantification}
\label{sec:sparse_grids}
% % % % % % % % % % % % % % % % % % % % 

% % % % % % % % % % % % % % % % % % % % 
\subsection{A surrogate-modeling approach to Uncertainty Quantification}
% % % % % % % % % % % % % % % % % % % % 

As already discussed in Section \ref{sec:model_problem}, the beam model 
depends on $N$ uncertain parameters, collected in the vector $\pp \in \Gamma$.
More precisely, we assume that each component $p_n$ is a uniform random variable that can take values in the range $\Gamma_n$ (we write $p_n\sim\mathcal U(\Gamma_n)$); 
we further assume that all random variables are independent, such that
the probability density function (pdf) of $\pp$ is simply the constant function  
$\rho(\pp) = \prod_{n=1}^N \frac{1}{b_n-a_n}$. 

Let us moreover denote by $\qoi \in \Rset$ the quantity of interest (QoI) or output of the beam equation (which we will call hereafter Full-Order Model, FOM), e.g., the displacement or the strain in a point of the beam.
$\qoi$ can then be seen as a $N$-variate function of the uncertain parameters, $\qoi = \qoi(\pp), \qoi: \Gamma \rightarrow \Rset$ (generalizations to vector-valued quantities of interest, i.e., $f: \Gamma \rightarrow \Rset^P$,
is straightforward; one such example is when we consider the entire displacement
field as QoI).

In this setup, we are interested in ``quantifying the uncertainties'' of the QoI due to the variability of $\pp$; 
to this end, we would like to compute statistical indices for $\qoi(\pp)$ such as its expected value and variance
\begin{align}
    & \Ev[\qoi] = \int_\Gamma \qoi(\pp) \rho(\pp) d\pp, \label{eq:exp-value} \\
    & \Var[\qoi] = \int_\Gamma (\qoi(\pp) - \Ev[\qoi])^2 \rho(\pp) d\pp = \Ev[\qoi^2] - \Ev[\qoi]^2, \nonumber
\end{align}
as well as higher order indices (such as kurtosis and skewness), and ideally its pdf. This task
is usually called UQ.

A successful approach to perform UQ is to build a so-called surrogate model for the QoI, following an offline/online paradigm. More precisely,
in a preliminary offline phase, a number of beam problems is solved, for certain judiciously selected combinations of values of $\pp$, and the corresponding values of $\qoi(\pp)$ stored; 
a so-called surrogate model is then constructed out of these values (by e.g. interpolation or
least-squares regression). During the subsequent online phase, quantities such as
those in Equation \eqref{eq:exp-value} are computed efficiently by evaluating the surrogate
model (cheap operation that essentially involves evaluating a polynomial expression) 
instead of repeatedly solving the beam problem. In the following, we construct a 
so-called \emph{sparse grids surrogate model}, but many other methods for building surrogate models are available in literature  (e.g. Polynomial Chaos, Reduced Basis, Gaussian Processes, Radial Basis Functions, Neural Networks, just to name a few - we refer e.g. to \cite{ghanem:UQbook} for an overview). In the context of timber engineering, surrogate models have also been employed in~\cite{klzwf_18}.

% % % % % % % % % % % % % % % % % % % % 
\subsection{Mathematical description of sparse grids} 
\label{sec:SG}
% % % % % % % % % % % % % % % % % % % % 

In this section, we quickly cover the basics of sparse grids, 
following closely the recent work \cite{piazzola.tamellini:SGK},
to which we refer the reader for more details.

The sparse grid surrogate model, 
which in the following will be denoted by $\SG_{\qoi}(\pp)$, 
can be informally described as an approximation of $\qoi(\pp)$, 
obtained as a linear combination of several ``small''
tensor interpolants of $\qoi$ over $\Gamma$, denoted $\qoi_{m(\ii)}(\pp)$ below,
each formed by a limited number of points. The underlying idea is the so-called
\emph{sparsification principle}, i.e., the intuition that, 
while none of these interpolants will be very accurate since they 
are all based on few points, 
by carefully combining many of them one can recover an overall good surrogate model. This comes 
at a much lower cost than what would be needed if one were to build naively a tensor
interpolant by covering the parameters space $\Gamma$ with a tensorial Cartesian
grid obtained by considering say $M$ values for each parameter. Indeed, such approach would involve a number of grid points exponential in the probabilistic dimension of the problems ($M^N$), i.e., it would be affected by the so-called \emph{curse of dimensionality}, which makes the tensor product technique unfeasible, even for even moderately small $N$. 

More precisely, the sparse grids surrogate model is expressed by means of the so-called \emph{combination technique} formula
\begin{equation}
\qoi(\pp) \approx S_{\qoi}(\pp) = \sum_{\ii \in \II} c_{\ii} \qoi_{m(\ii)}(\pp),
\quad
c_{\ii} = \sum_{\substack{\jj \in \{0,1\}^N \\ \ii+\jj \in \II}} (-1)^{|\jj|}, 
\label{eq:sparsegrid-combitec}
\end{equation}  
where:
\begin{itemize}
    \item $\ii \in \Nset^N$ is a multi-index, i.e., a vector of $N$ integer numbers; 
    a tensor interpolant $\qoi_{m(\ii)}(\pp)$ will be associated to each $\ii$ 
    in the set $\II$ (more on this later), and each entry $i_k$ of $\ii$ 
    denotes the \emph{level} of approximation of $\qoi_{m(\ii)}(\pp)$
    along each parameter $p_k, k=1,\ldots N$;
    \item $m(\cdot)$ is an increasing function (``level-to-knots function''), 
        such as $m(k)=k$ or $m(k)=2^k$;
    \item $m(\ii)$ is the vector obtained applying $m(\cdot)$ to each component of $\ii$,
        i.e., $m(\ii)=[m(i_1),m(i_2),\ldots]$; 
    \item $\qoi_{m(\ii)}(\pp)$ is a tensor interpolant, built over a Cartesian 
        grid on $\Gamma$ with $m(i_1) \times m(i_2) \times \cdots$ points;
        more details on the construction and evaluation of $\qoi_{m(\ii)}(\pp)$
        are reported in \ref{appendix:sparse_grids_formulas}.
    \item $c_{\ii}$ are the so-called \emph{combination technique coefficients}. 
        Note that some $c_{\ii}$ might be null, in which case $\qoi_{m(\ii)}(\pp)$ is not
        part of the final approximation;
    \item $\II$ is a multi-index set, $\II \subset \Nset^N$, that specifies which tensor     
        interpolants are candidates to enter in the sparse grid construction.
        It should be chosen according to the sparsification principle mentioned above,
        and in particular it should refrain from containing indices $\ii$
        whose entries are all large numbers
        (the cost of building the associated interpolant $\qoi_{m(\ii)}(\pp)$
        would be too large). Instead, whenever one entry (or a few entries) of $\ii$ are 
        large, the others should be kept as small as possible.
        Moreover, for technical reasons it is required that $\II$ is \emph{downward-closed}, i.e., 
        if $\ii \in \II$ then all is ``precedent'' 
        neighbors are also in $\II$\footnote{
        Upon denoting with $\ee_k$ the $k$-th versor of
        $\Nset^N$, i.e., the vector with all zeros expect the $k$-th component, that is equal 
        to 1, $\II$ is downward-closed if:  
        \[
            \ii \in \II \Rightarrow \ii - \ee_k \in \II, \forall k = 1,\ldots,N    
        \]}.
\end{itemize}
The set of points where $\qoi(\pp)$ is evaluated (i.e., the 
union of all the points needed to assemble each $\qoi_{m(\pp)}(\pp)$ ) is called
\emph{sparse grid}, and will be denoted by $\GGG$. Its cardinality will be denoted by $G$.

\begin{figure*}
    \centering
     \subfigure[]{\includegraphics[width=0.48\textwidth]{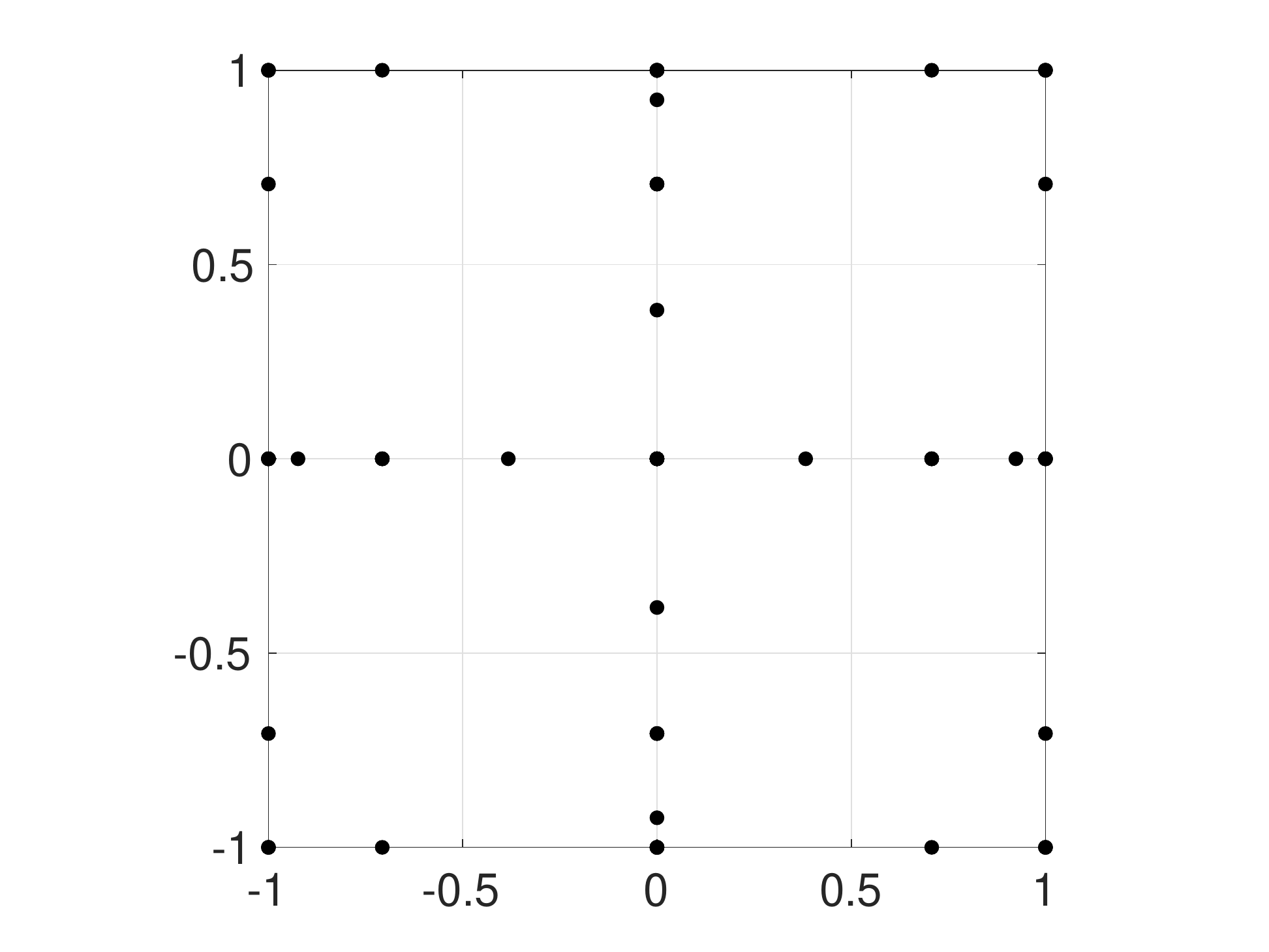}
        \label{fig:sg2D_w3}}
    \subfigure[]{\includegraphics[width=0.48\textwidth]{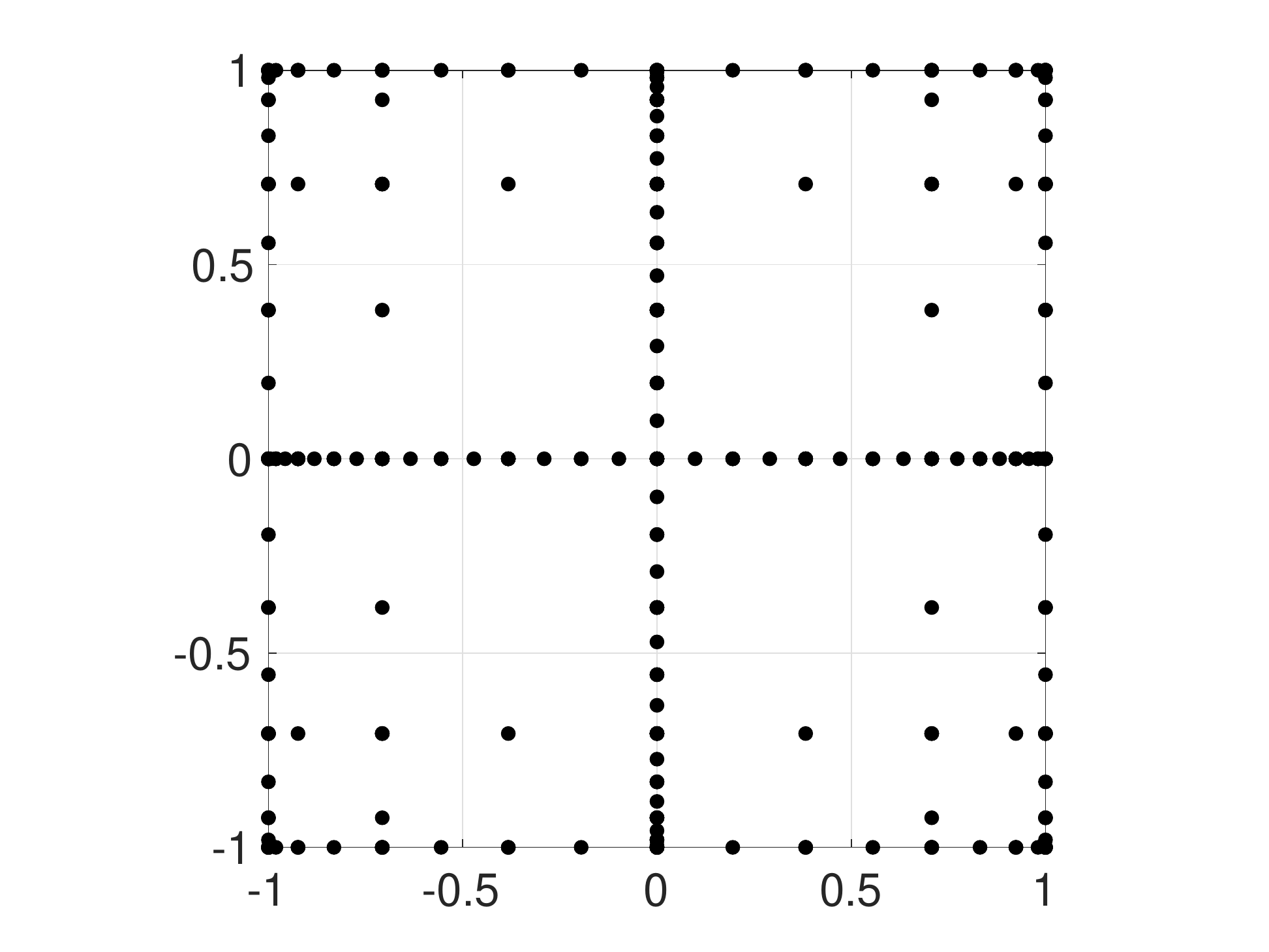}
        \label{fig:sg2D_w5}}
 
     \caption{
     Examples of sparse grids of level $w=3$ (Figure~\ref{fig:sg2D_w3}) and level $w=5$ (Figure~\ref{fig:sg2D_w5}) in $N=2$ dimensions. We have used Clenshaw-Curtis points~\eqref{eq:points}, function $m(\cdot)$ as in~\eqref{eq:m} and multi-index set $\emph{I}$ as in~\eqref{eq:I}.}
\end{figure*}

Equation \eqref{eq:sparsegrid-combitec} becomes operative the moment we specify the three basic ``ingredients'' of the sparse grid construction, namely, the set $\II$, the function $m(\cdot)$ and the knots used to construct each tensor interpolant $\qoi_{m(\ii)}(\pp)$. A lot of literature deals with criteria and algorithms to optimally choose these three components. In Section~\ref{sec:numerical_experiments} we detail the choices we have adopted in this work. For examples of sparse grids in $N=2$ dimensions we refer to Figures~\ref{fig:sg2D_w3} and~\ref{fig:sg2D_w5}.

% % % % % % % % % % % % % % % % % % % % 
\subsection{Sparse grids for UQ} 
\label{sec:sparse_grids_for_UQ}
% % % % % % % % % % % % % % % % % % % % 

In this section we provide an overview on how sparse grids can be used for UQ of a quantity of interest $\qoi$. In our work, we have used the implementation of sparse grids provided in the Sparse Grids Matlab Kit, 
see \cite{piazzola.tamellini:SGK}, which can be used essentially ``off-the-shelf'' and
renders these operations rather straightforward (this can be 
appreciated also by taking a look at the Listings reported in Section~\ref{sec:numerical_experiments}).

\begin{description}
    \item[Expected value.] 
    The univariate interpolation points used as basic blocks of the sparse grid construction
    always come associated with quadrature weights. For example, the weights corresponding to Clenshaw--Curtis points (used in our numerical experiments, 
    see Section \ref{sec:numerical_experiments}) can be computed by Fast Fourier transform \cite{trefethen:comparison}.
    Recalling that expected values are just weighted integrals over $\Gamma$ 
    (cf. Equation \eqref{eq:exp-value}),
    a \emph{sparse grid quadrature} $\mathcal{Q}[f]$ can be derived
    (mimicking the steps that would lead to Equation \ref{eq:sparsegrid-combitec}),
    which in practice simply amounts to taking weighted
    sums 
    %leading eventually to the approximation of the expected value as a linear combination 
    of the evaluations of $\qoi$ over the sparse grid points $\qq \in \GGG$. 
    The weights $\alpha_{\qq}$ depend on quadrature weights of the interpolation points
    and on the combination technique coefficients $c_{\ii}$ 
    (see \cite{piazzola.tamellini:SGK} for details):
    \[
    \Ev [\qoi] = \int_{\Gamma} \qoi(\pp) \rho(\pp) d\pp 
    \approx \sum_{\qq \in \GGG} \alpha_{\qq} \qoi(\qq) = 
    \mathcal{Q}[f].
    \]
    See Listing \ref{ls:QoI} for software calls.
    
    \item[Variance and higher order indices.] 
    Simply employ the fact 
    already recalled in Equation \eqref{eq:exp-value} that
    $\Var[f] = \Ev[f^2] - \Ev[f]^2$,
    and approximate both terms by sparse grids quadrature as explained above.
    Similar formulas exists for higher moments such as skewness (connected to $\Ev[f^3]$)
    and kurtosis (connected to $\Ev[f^4]$). 
    
    \item[Global sensitivity analysis by Sobol indices.] 
    Sobol indices \cite{sudret:sobol,archer.saltelli.sobol:anova}
    are quantities that
    assess the contribution of each uncertain parameter to the total variance of a 
    quantity of interest; the underlying mathematical machinery is a decomposition of 
    the variance of $\qoi$ similar to the ANOVA decomposition. In particular, 
    the \emph{principal} Sobol index $S_i^P$ 
    quantifies the impact of each uncertain parameter $p_i$ alone,
    whereas the \emph{total} Sobol index $S_i^T$ 
    quantifies the impact of each uncertain parameter 
    alone and in mixed effect with any other uncertain parameter. 
    Principal and Sobol indices can be obtained by post-processing of the sparse grid
    surrogate model $\SG_{\qoi}(\pp)$, see \cite{feal:compgeo} for details. See Listing \ref{ls:sobol} for implementation details.

    \item[Probability density function.] 
    An approximation of the pdf of $\qoi$ can be obtained by generating
    sufficiently many samples of the uncertain parameters $\pp_i$ according to their 
    pdf, evaluating $\qoi$ for each of them, and then resorting to binning algorithms
    to generate histograms of such values, or using functions such as 
    kernel density estimates \cite{rosenblatt:kde}. 
    This process is significantly sped up by replacing 
    the values $\qoi(\pp_i)$ with their approximate counterparts $\SG_{\qoi}(\pp_i)$,
    \cite{sagiv:pdfconv}.
    To this end we remark that evaluating  $\SG_{\qoi}(\pp_i)$ is essentially real-time
    (one only needs to evaluate a few polynomial interpolants) whereas evaluating 
    $\qoi(\pp_i)$ requires solving a PDE (beam problem).
    See Listing \ref{ls:pdf} for implementation details.
\end{description}

% % % % % % % % % % % % % % % % % % % % 
\section{Numerical experiments}
\label{sec:numerical_experiments}
% % % % % % % % % % % % % % % % % % % % 

\begin{figure}
    \centering
    \includegraphics[width=0.48\textwidth]{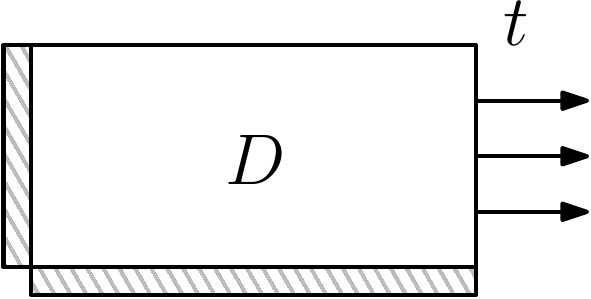}
    \caption{Traction model. The external load is $\tt=(1 \kilo\newton/\meter,0)^T$, homogeneous Neumann boundary conditions are imposed at the right and top part of the boundary, whereas the horizontal (vertical, respectively) displacement is imposed zero at the left (bottom, respectively) part of the boundary.}
    \label{fig:domain}
\end{figure}

All the numerical tests deal with the traction problem (see Figure~\ref{fig:domain}), namely we take the external load $\tt=(10^3 \kilo\newton/\meter,0)^T$ and we impose homogeneous Dirichlet boundary conditions on $\Sigma_s=[0,L]\times\{0\}\cup\{0\}\times[0,H]$ and homogeneous Neumann boundary conditions on $\partial D\setminus\Sigma_s$. 

\bigskip

We detail now the choices adopted for the sparse grid construction illustrated in Section~\ref{sec:SG}.
\begin{itemize}
    \item As knots, we use the Clenshaw--Curtis points  
    points), which are well suited for uncertain parameters with uniform pdf
    \footnote{Note that equispaced points are in general \emph{not} a good choice,
    due to the well-known Runge's phenomenon.}.
    A set of $K$ points in $[-1,1]$ can be computed as follows:
    \begin{equation}
        \label{eq:points}
        x_K^{(j)} =\cos\left(\frac{(j-1) \pi}{K-1}\right), \quad 1\leq j \leq K,    
    \end{equation}
    and then if needed linearly transformed to any generic interval $[a,b]$.   
    \item As function $m(\cdot)$, we use the following
    \begin{equation}
    \label{eq:m}
        m(k) = 
    \begin{cases}
    1, & k=1, \\
    2^{k-1}+1, & k>1 ,  
    \end{cases}    
    \end{equation}
    which yields the doubling of the number of interpolation points, when moving from interpolation level $k$ to $k+1$. Note that this is choice is particularly useful since it renders tensor interpolants built with Clenshaw--Curtis points \emph{nested}, i.e.,
    the set of points needed to build $\qoi_{m(\jj)}(\pp)$ is contained in the
    set needed to build $\qoi_{m(\ii)}(\pp)$ if the multi-indices $\ii=(i_1,\ldots,i_N)$ and $\jj=(j_1,\ldots,j_N)$ fulfill  $j_k \leq i_k$ for all $k=1,\ldots,N$. 
    This property is clearly beneficial if one wants to refine a sparse grid already 
    computed by adding further computations.
    \item As set $\II$ we use the classical choice
    \begin{equation}
        \label{eq:I}
        \II = \left\{ \ii \in \Nset^N : t(\ii) \leq w\right\},    
    \end{equation}
    where $t\colon\Nset^N\rightarrow\Nset$ is given by $t(\ii)\coloneqq \sum_{n=1}^N i_n$ and $w\in \Nset$ is an integer number that controls the accuracy of the sparse grid
    (the larger $w$, the more points in the sparse grid). It is easy to see that
    it enforces a basic version of the sparsification principle;
    more sophisticated choices, such as anisotropic sets or adaptive algorithms
    for the selection of $\II$ are discussed e.g. in
    \cite{nobile.tempone.eal:aniso,chkifa:adaptive-interp,nobile.eal:optimal-sparse-grids}
    and \cite{gerstner.griebel:adaptive,nobile.eal:adaptive-lognormal}, respectively.
\end{itemize}
These three choices combined generate a sparse grid, which is commonly named in the literature as \emph{Smolyak grid}. 

\bigskip

In the following we discuss two numerical examples with load and boundary conditions as in Figure \ref{fig:domain}.
\begin{itemize}
\item     
In Section \ref{sec:one_knot} we use three uncertain parameters to model the presence of one knot inside the unit square domain.
The elasticity modulus fulfills the simplified assumption of being $y$-independent. We consider in this example two QoI, namely
(i) the entire displacement field (ii) the horizontal displacement at the bottom-right corner of the domain.
The outcomes are: (i) the numerical study of the approximation error of the first QoI; (ii) the construction of the surrogate for the second QoI and the numerical study of its accuracy; (iii) pdf and Sobol indices of the second QoI.

\item 
In Section \ref{sec:two_knot} we consider a rectangular domain and we use seven uncertain parameters to model the presence of two knots. In contrast to Section \ref{sec:one_knot}, here the elasticity modulus varies along both the horizontal and the vertical directions.
Differently from before, in this example we consider only one
QoI, i.e., the horizontal displacement at the bottom-right corner of the domain, and we compare two surrogates computed by means of Smolyak sparse grids and a-posteriori adaptive sparse grids,
i.e., different strategies to compute the set $\II$.
\end{itemize}

% % % % % % % % % % % % % % % % % % % % 
\subsection{One-knot example}
\label{sec:one_knot}
% % % % % % % % % % % % % % % % % % % % 

In the first numerical example we take $L = H = 1 \meter$ and we choose the stochastic elasticity modulus~\eqref{eq:E} depending on the uncertain vector $\pp=(p_1,p_2,p_3)$ with length $N=3$. 
More in details, we take $E_0=10^4 \mega\pascal$ and 
\begin{equation}
    \label{eq:alpha_1knot}
    \alpha(x,y,\pp)
    = p_1 - \gamma \exp\left(-\frac{(x-p_2)^2}{2p_3^2}\right),
\end{equation}
with $p_1\sim\mathcal U(0.5,1.5)$, $p_2\sim\mathcal U(0.25,0.75)$, $p_3\sim\mathcal U(0.1,0.2)$ and $\gamma=0.4$. 

\begin{figure*}
    \centering
    \subfigure[]
    {\includegraphics[width=0.48\textwidth]{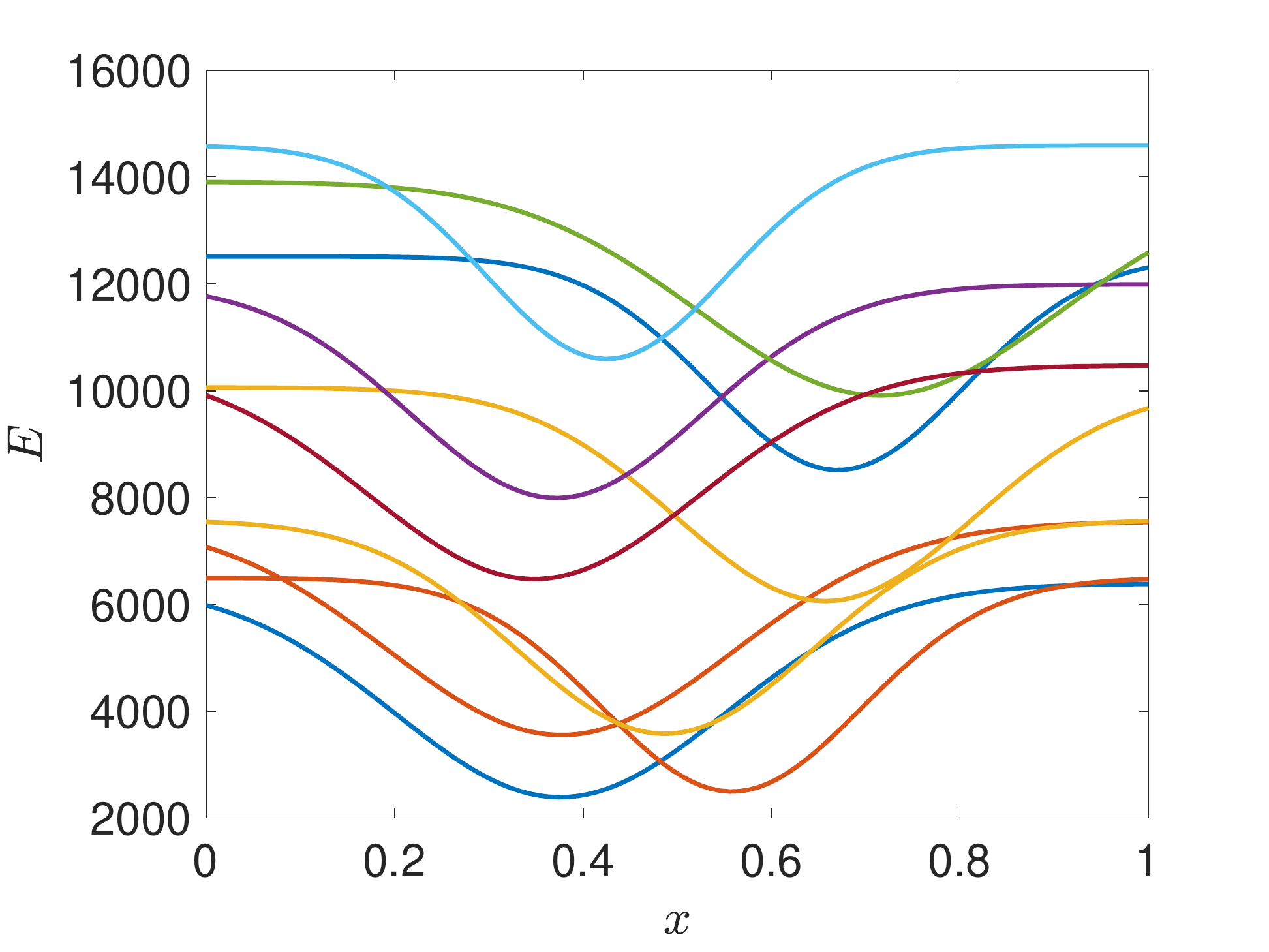}
        \label{fig:alpha_1knot_section}}
    \subfigure[]
    {\includegraphics[width=0.48\textwidth]{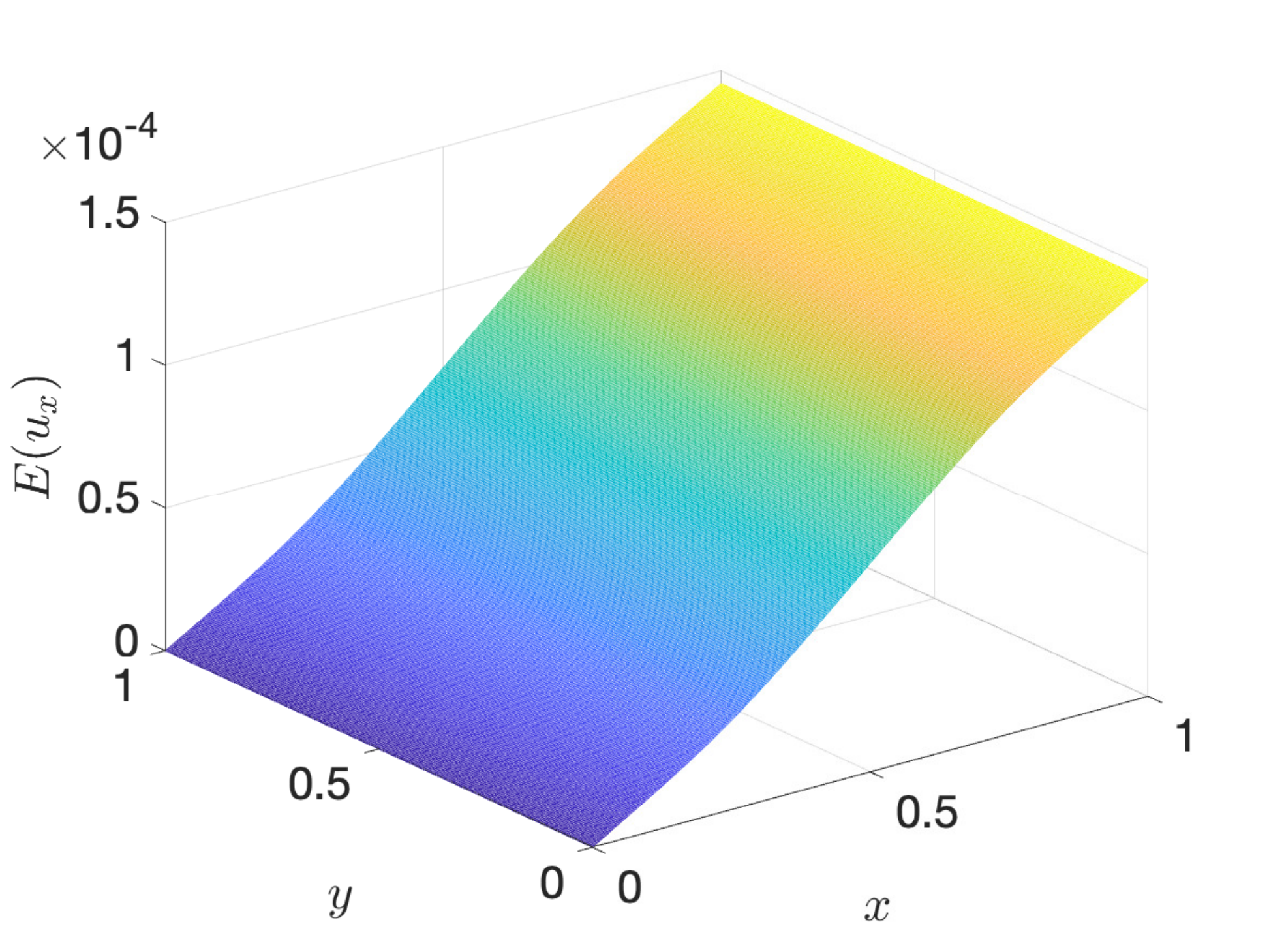}
        \label{fig:mean_ref_1knot}}
    
    \subfigure[]
    {\includegraphics[width=0.48\textwidth]{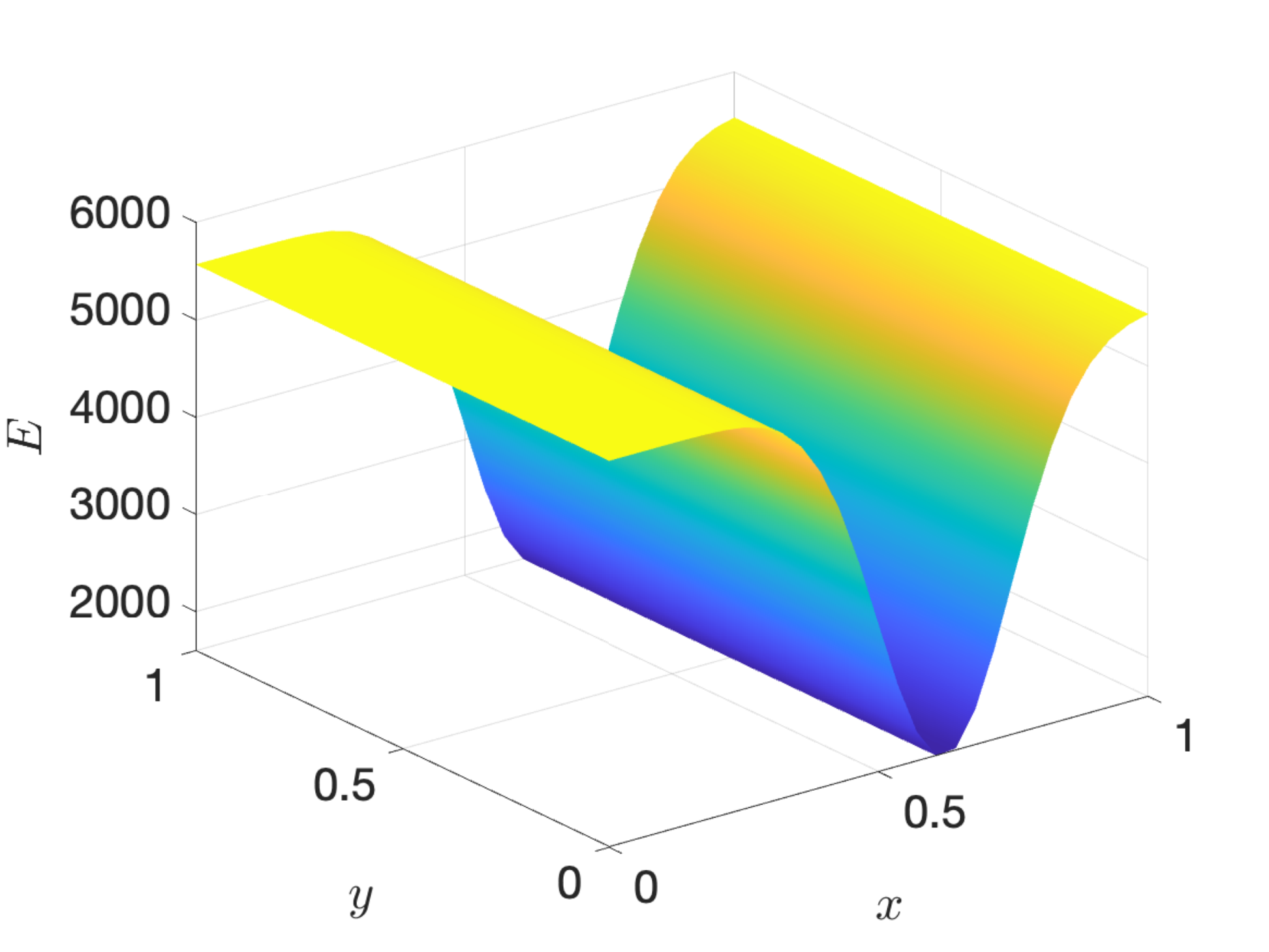}
        \label{fig:alpha_1knot_sample1}}
    \subfigure[]
    {\includegraphics[width=0.48\textwidth]{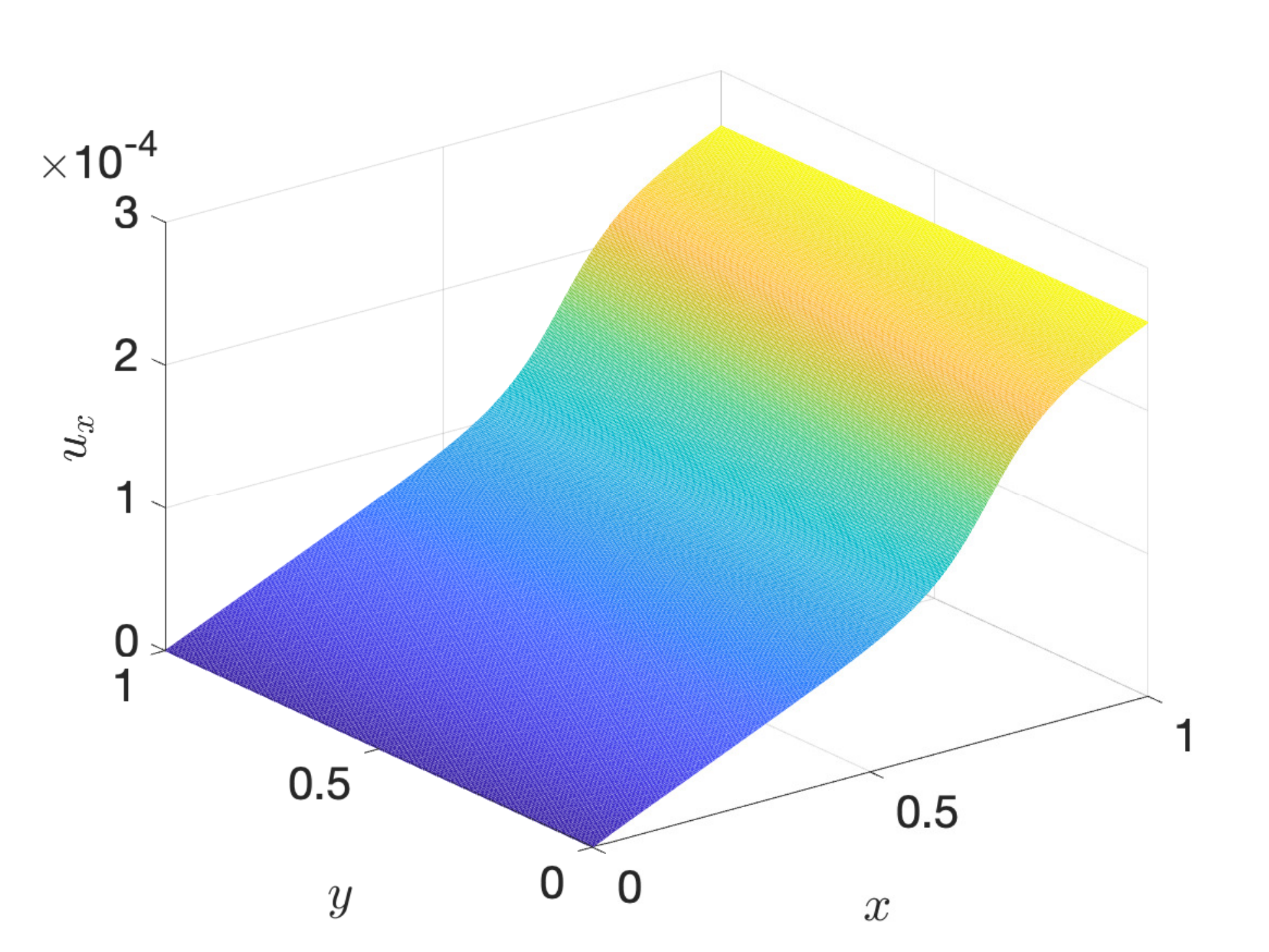}
        \label{fig:sol_1knot_sample1}}
    
    \subfigure[]
    {\includegraphics[width=0.48\textwidth]{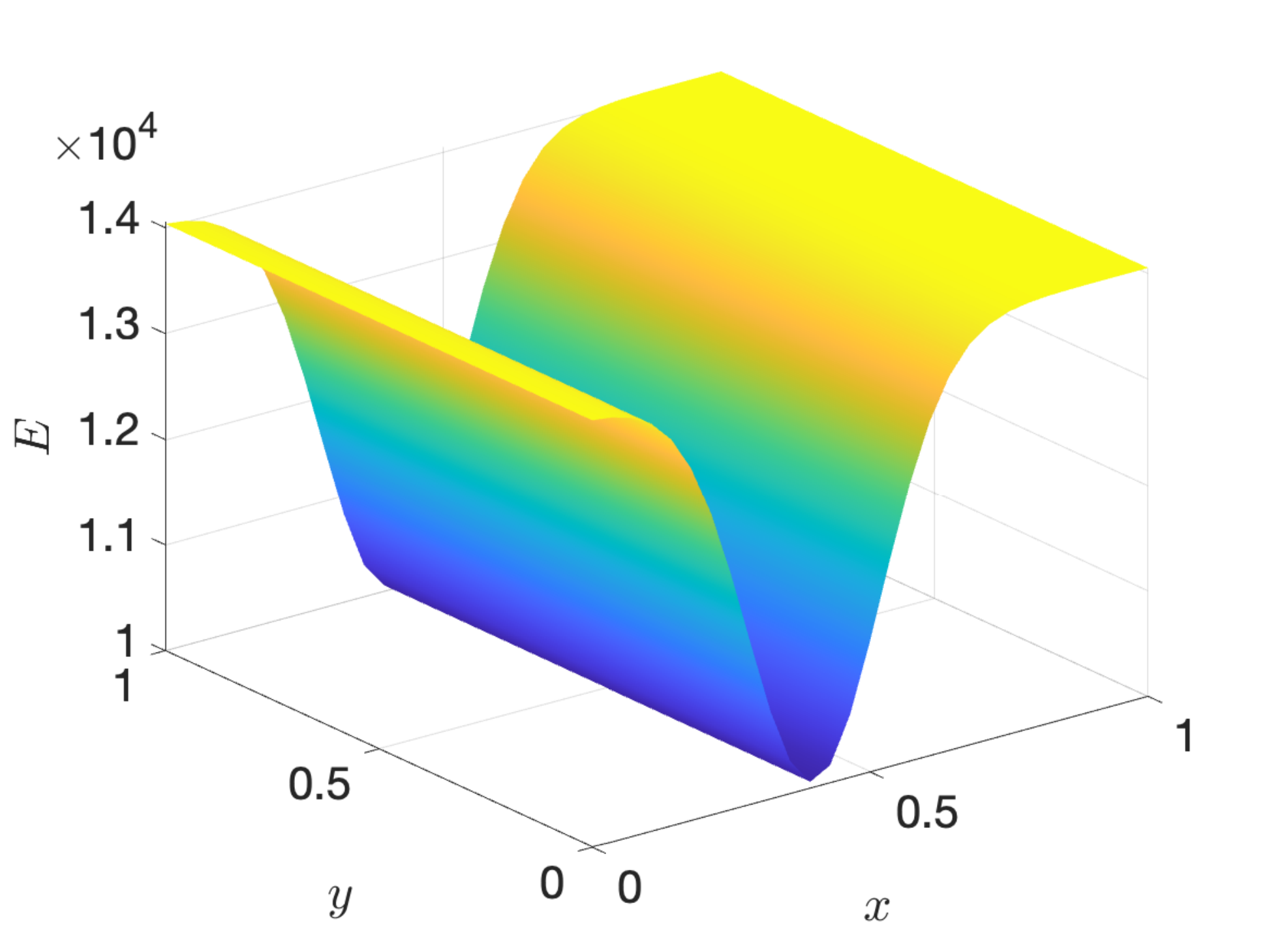}
        \label{fig:alpha_1knot_sample2}} 
    \subfigure[]
    {\includegraphics[width=0.48\textwidth]{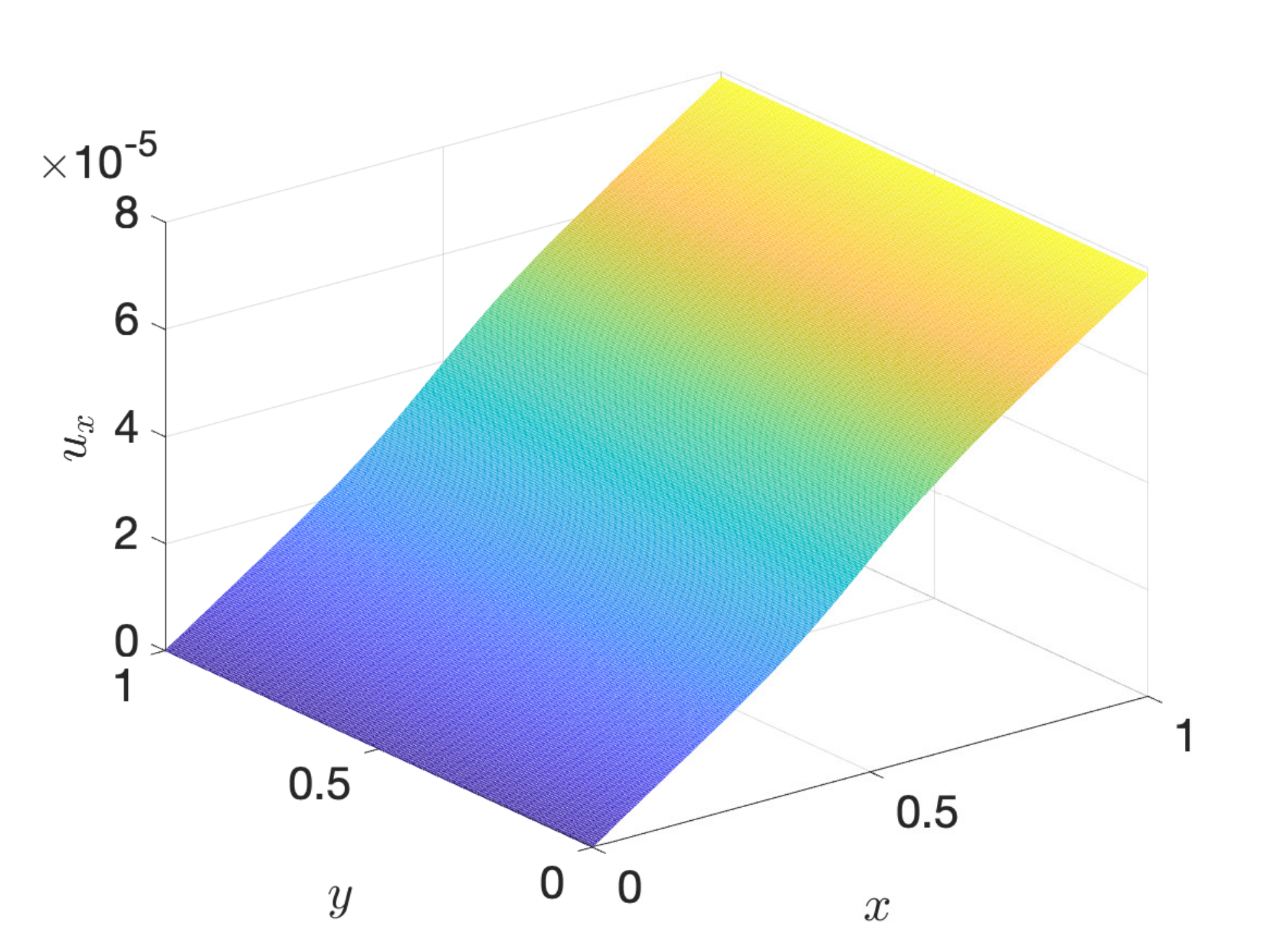}
        \label{fig:sol_1knot_sample2}}
    
    \caption{
        \ref{fig:alpha_1knot_section} Ten samples of $E(x,\bar y,\pp)$, for fixed $\bar y$;
        \ref{fig:mean_ref_1knot} Reference solution for $\Ev[u_x]$ (Smolyak sparse grid surrogate of level $w=7$);
        \ref{fig:alpha_1knot_sample1}, \ref{fig:alpha_1knot_sample2} Plot of two samples of $E(x,y,\pp)$ as in~\eqref{eq:alpha_1knot} for $\pp=(5567,0.62277,0.12425)$ and $\pp=(14052,0.39997,0.11967)$, respectively;
        \ref{fig:sol_1knot_sample1}, \ref{fig:sol_1knot_sample2} Plot of the two corresponding solutions $u_x(x,y,\pp)$~\eqref{eq:pde} ($\pp=(5567,0.62277,0.12425)$ and $\pp=(14052,0.39997,0.11967)$, respectively) computed via IGA (Section~\ref{sec:IGA}).}
    
\end{figure*}

Note that $\alpha(\cdot,\pp)\colon D\rightarrow\mathbb R_+$ varies in the horizontal direction $x$, only, whereas it is constant in the vertical direction $y$. Therefore, the displacement along the vertical direction $u_y$ is zero. 
This choice of $\alpha$ aims at modeling the presence of one knot along the beam. Following this interpretation, $p_2$ represents the (variable) center of the knot and $p_3$ represents its (variable) width; finally, $p_1E_0$ is the (variable) nominal value of the Young modulus away from the knot. Note that the ranges of $p_2$ and $p_3$ are chosen so that the knot is well-contained inside the beam.
We refer to Figure~\ref{fig:alpha_1knot_section}, depicting a set of ten samples of $E$ plotted versus the horizontal variable $x\in [0,1]$, and Figure~\ref{fig:alpha_1knot_sample1}, Figure~\ref{fig:alpha_1knot_sample2}, where two samples of $E$ are plotted versus $(x,y)\in D$. 

The IGA approximation of the corresponding solutions of problem~\eqref{eq:pde} are shown in Figure~\ref{fig:sol_1knot_sample1} and Figure~\ref{fig:sol_1knot_sample2}. For this numerical experiment, the IGA parameters are set as $\degr=\degq=4$ and $\ncp=\mcp=32$, leading to a negligible error in the space variables. 

With the above-mentioned choices, the sparse grid of level $w$ can be generated by running the very simple Matlab code in Listing~\ref{ls:sparse_grid}.

% % % % % % % % % % % % % % % % % % % % 
\subsubsection{QoI 1: displacement field}
% % % % % % % % % % % % % % % % % % % % 

The surrogate for the solution $\uu$ or a QoI of $\uu$ can then be computed by running the Matlab code in Listing~\ref{ls:QoI}. 
The brevity and simplicity of these listings testify how little extra work is needed to
interface the beam solver to the UQ software, and thus how easy it is to perform a UQ analysis.

\begin{lstlisting}[caption = {Matlab code to create the Smolyak sparse grid.}, label={ls:sparse_grid}]
% number of parameters
N=3; 
% knots for p1, p2 and p3
knots_p1=@(n) knots_CC(n,0.5,1.5,'prob'); 
knots_p2=@(n) knots_CC(n,0.25,0.75,'prob'); 
knots_p3=@(n) knots_CC(n,0.1,0.2,'prob'); 
knots = {knots_p1,knots_p2,knots_p3};
% functions m and t
[lev2knots,idxset]=define_functions_for_rule('SM',N)
% level
w = 1; 
% sparse grid
S = smolyak_grid(N,w,knots,lev2knots,idxset,[]); 
% creates the "uniqued" list of points
Sr = reduce_sparse_grid(S) 
\end{lstlisting}

\begin{lstlisting}[caption = {Matlab code to compute the surrogate for the QoI.}, label={ls:QoI}]
% wrap the beam IGA solver into an @-function
f = @(y) solve_PDE(y); 
% or further have the solver just return the QoI
% f = @(y) QoI(y);
Eu = quadrature_on_sparse_grid(@(y)f(y),S,Sr);
\end{lstlisting}

\begin{figure}
    \centering
    \includegraphics[width=0.48\textwidth]{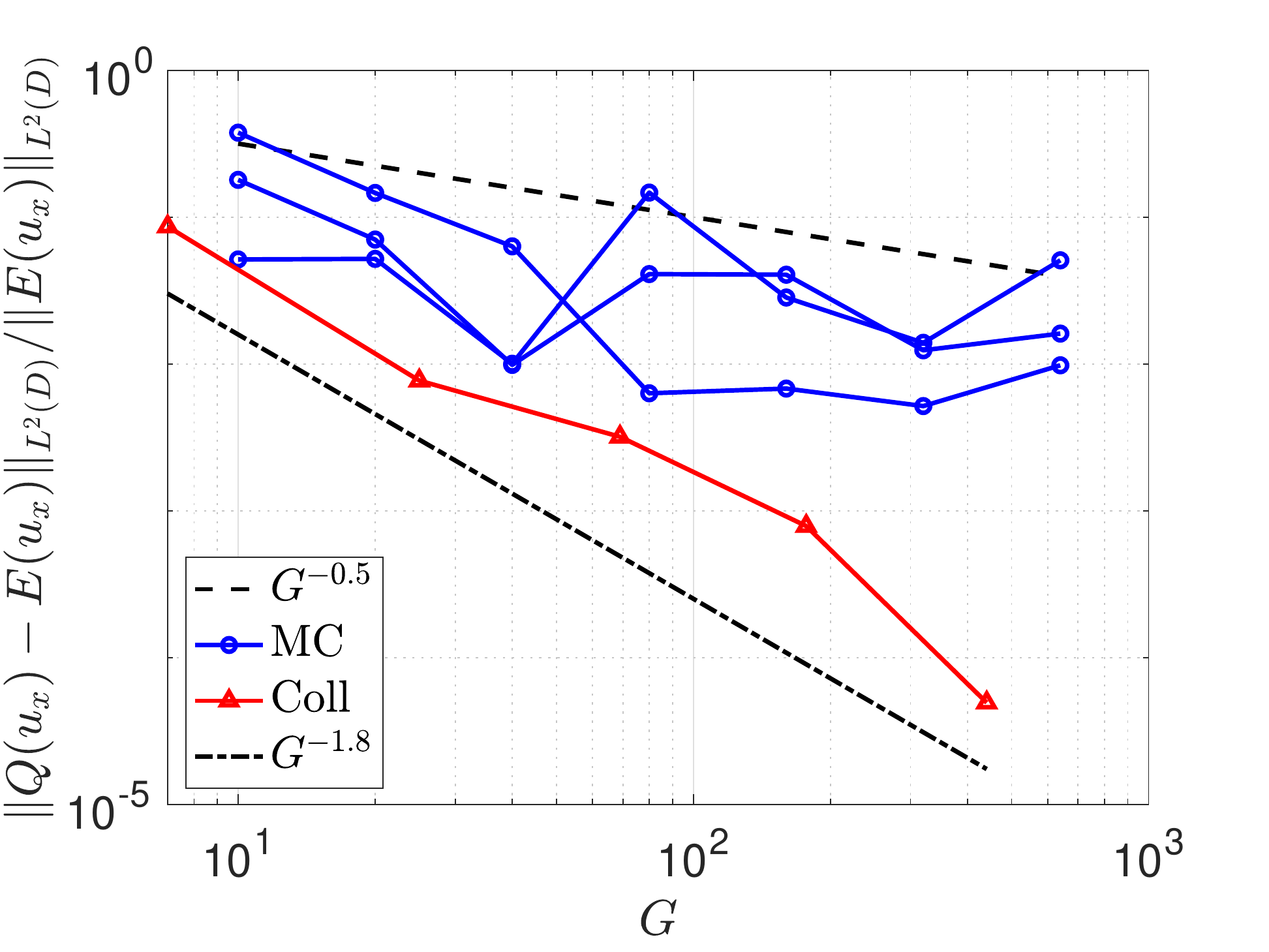}
     \caption{
     Error convergence of the expected value of the horizontal displacement field in the $L^2(D)$-norm.
     The quantities are plotted versus the number of PDE solves (the cardinality $G$ of the sparse grid for the Collocation method, the number of samples for the Monte Carlo method).}
    \label{fig:error_mean}
\end{figure}

We are interested in particular in assessing the quality 
of the approximation of the expected value of the 
horizontal displacement field, i.e., of $\Ev[u_x]$, 
which we approximately obtained employing the Smolyak sparse grid of level $w=7$ (see Figure~\ref{fig:mean_ref_1knot}). 
Coarser approximations to the expectation of the same quantity are then computed by Smolyak sparse grids of lower levels $w=1,\ldots,5$. Their relative error with respect to the
reference solution, measured in the $L^2(D)$-norm, is plotted in Figure~\ref{fig:error_mean}: the horizontal axis reports the cardinality of the employed sparse grid. For the sake of comparison, three instances of convergence of the Monte Carlo method are also depicted. 
When the sparse grid method is employed we observe an algebraic decay of the error with estimated slope -1.8, as opposed to the usual Monte Carlo decay rate $-1/2$ (i.e., the inverse of the square root of the number of Monte Carlo samples). We underline the effectiveness of the sparse grid approach, which delivers more accurate results with many less samples points as the plain Monte Carlo method.

\begin{lstlisting}[caption = {Matlab code to plot the sparse grid surrogate $\mathcal S_\qoi$.}, label={ls:plot}]
% define range of the parameters
aa = [0.5, 0.25, 0.1];
bb = [1.5, 0.75, 0.2];
domain = [aa; bb];
f_values = evaluate_on_sparse_grids(f,Sr); 
plot_sparse_grids_interpolant(S,Sr,domain,f_values,'with_f_values');
\end{lstlisting}

% % % % % % % % % % % % % % % % % % % % 
\subsubsection{QoI 2: horizontal displacement at the bottom-right corner}
% % % % % % % % % % % % % % % % % % % % 

\begin{figure*}
    \centering
    \subfigure[]{\includegraphics[width=0.48\textwidth]{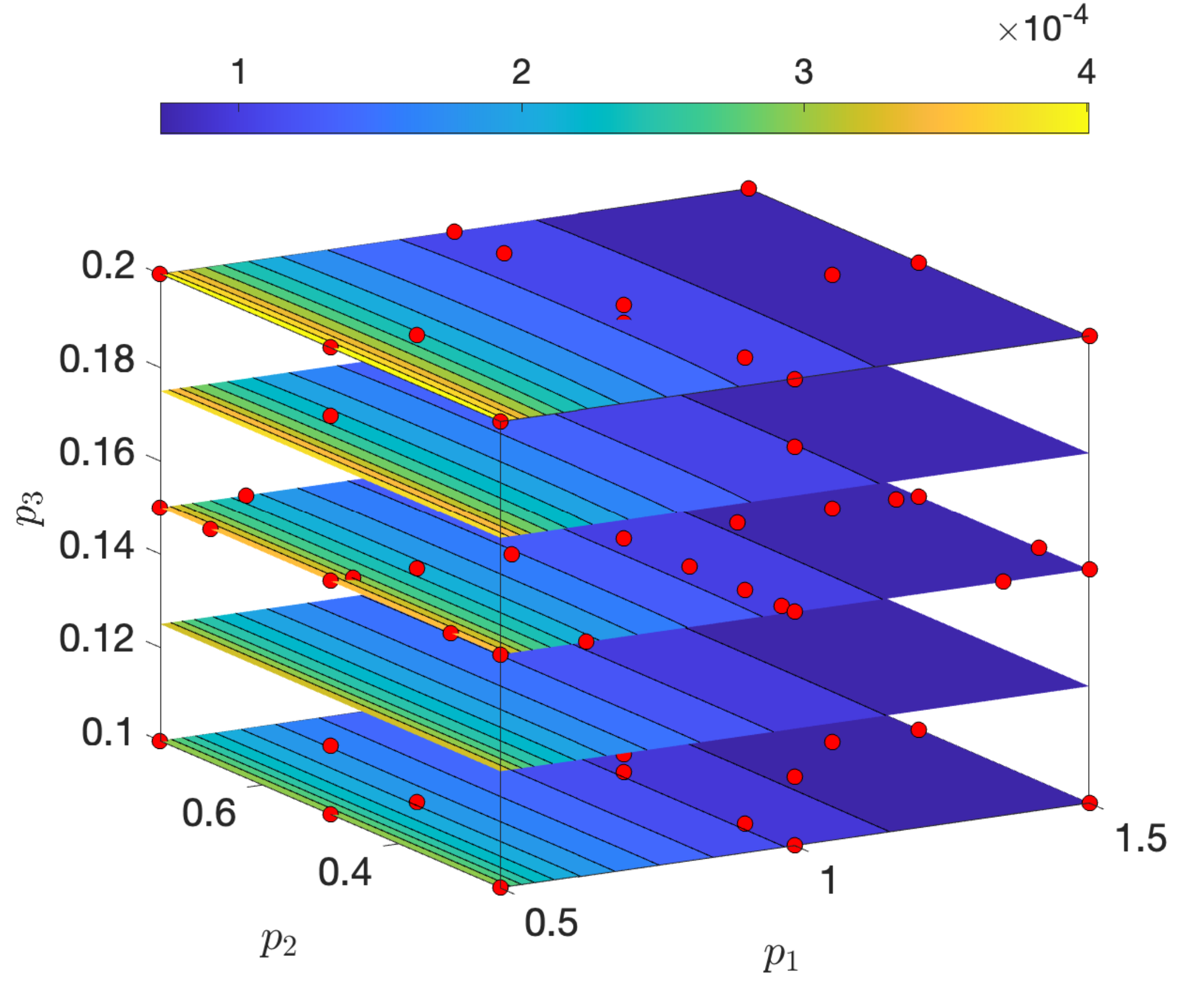}
        \label{fig:surface_1knot}}
    \subfigure[]{\includegraphics[width=0.48\textwidth]{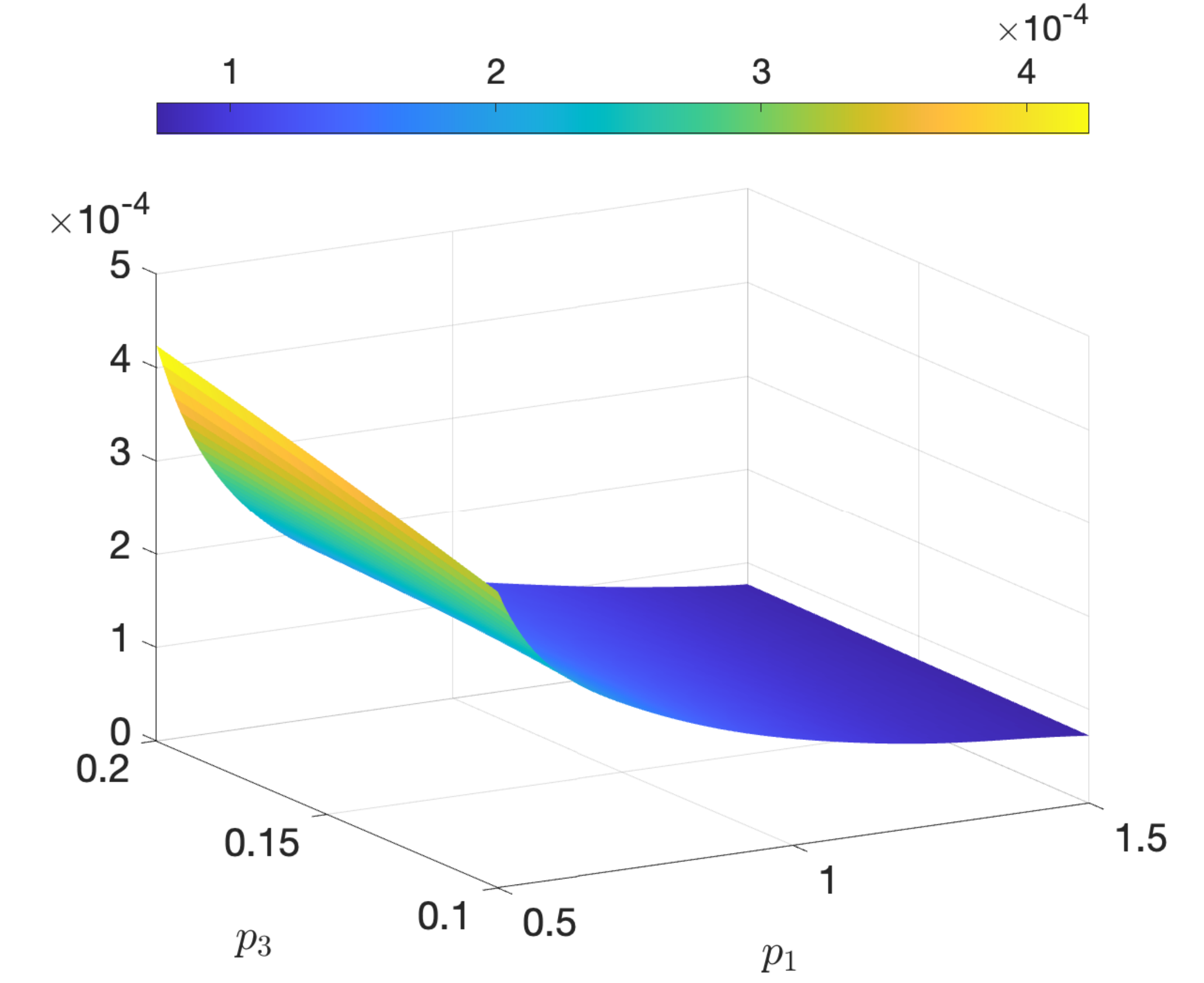}
        \label{fig:surface_p2fixed_1knot}}    
        
     \caption{
     \ref{fig:surface_1knot} Three-dimensional plot of the surrogate QoI computed on the sparse grid with level $w=3$;
     \ref{fig:surface_p2fixed_1knot} Plot of the surrogate QoI versus $p_1,p_3$ and for fixed $p_2=0.25$.}
     \label{fig:surface}
\end{figure*}

Let us consider now the real-valued QoI being the evaluation of the horizontal displacement at the bottom-right corner of the beam, namely $\qoi(\pp)=u_x(1,0,\pp)$.
The surrogate of the QoI can be easily computed (see Listing \ref{ls:QoI}) and plotted: see Listing~\ref{ls:plot} and Figure~\ref{fig:surface_1knot}, where level $w=3$ is considered. In Figure~\ref{fig:surface_p2fixed_1knot} we display a section of the three-dimensional plot in Figure~\ref{fig:surface_1knot} obtained for the fixed value $p_2=0.25$. The plot shows that the variability of the QoI with respect to the parameters $p_2$ and $p_3$ is very limited. This observation will be confirmed later on, by means of the Sobol indices.

\begin{figure*}
    \centering
     \subfigure[]{\includegraphics[width=0.48\textwidth]{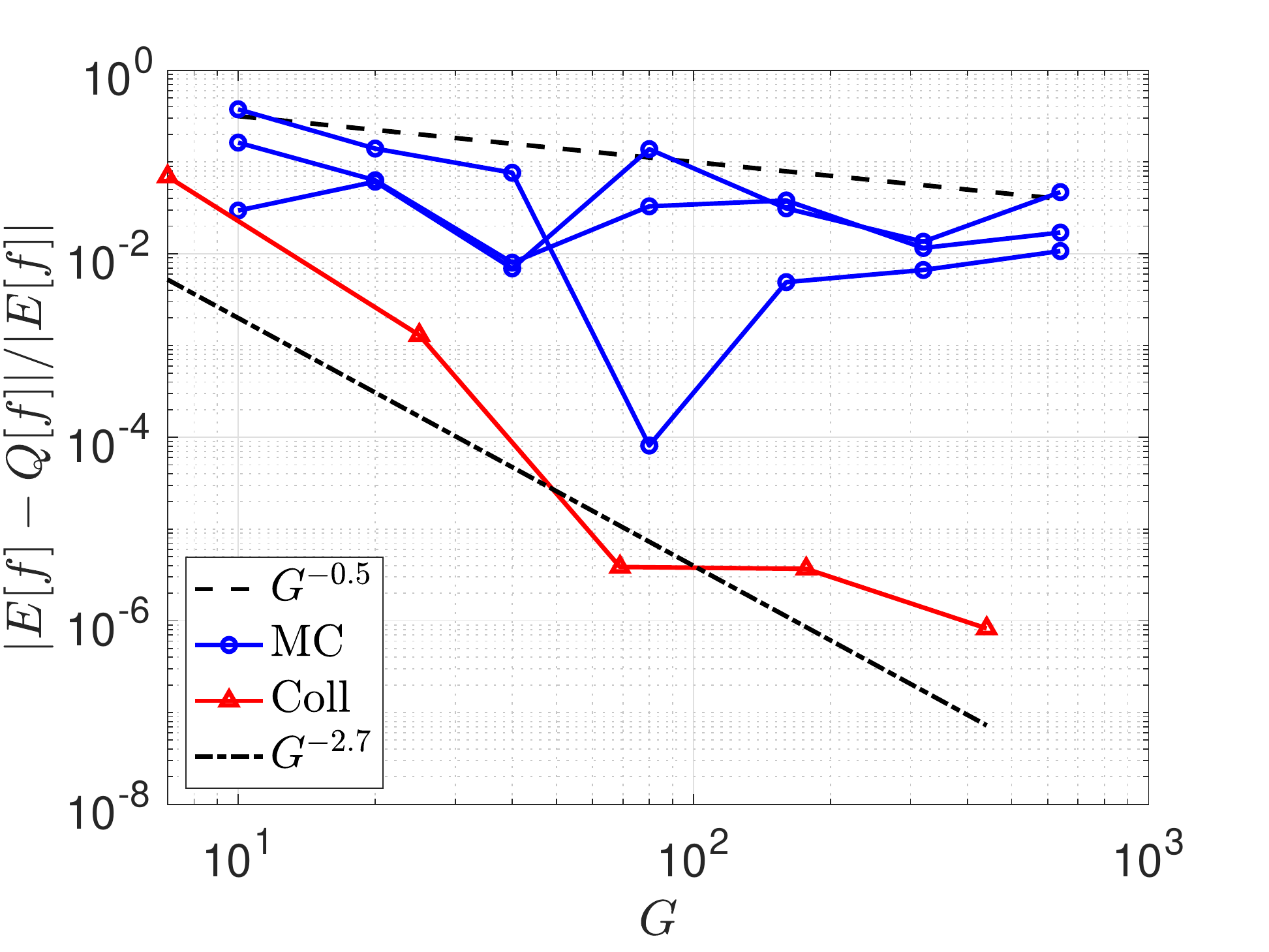}
        \label{fig:error2_mean}}
    \subfigure[]{\includegraphics[width=0.48\textwidth]{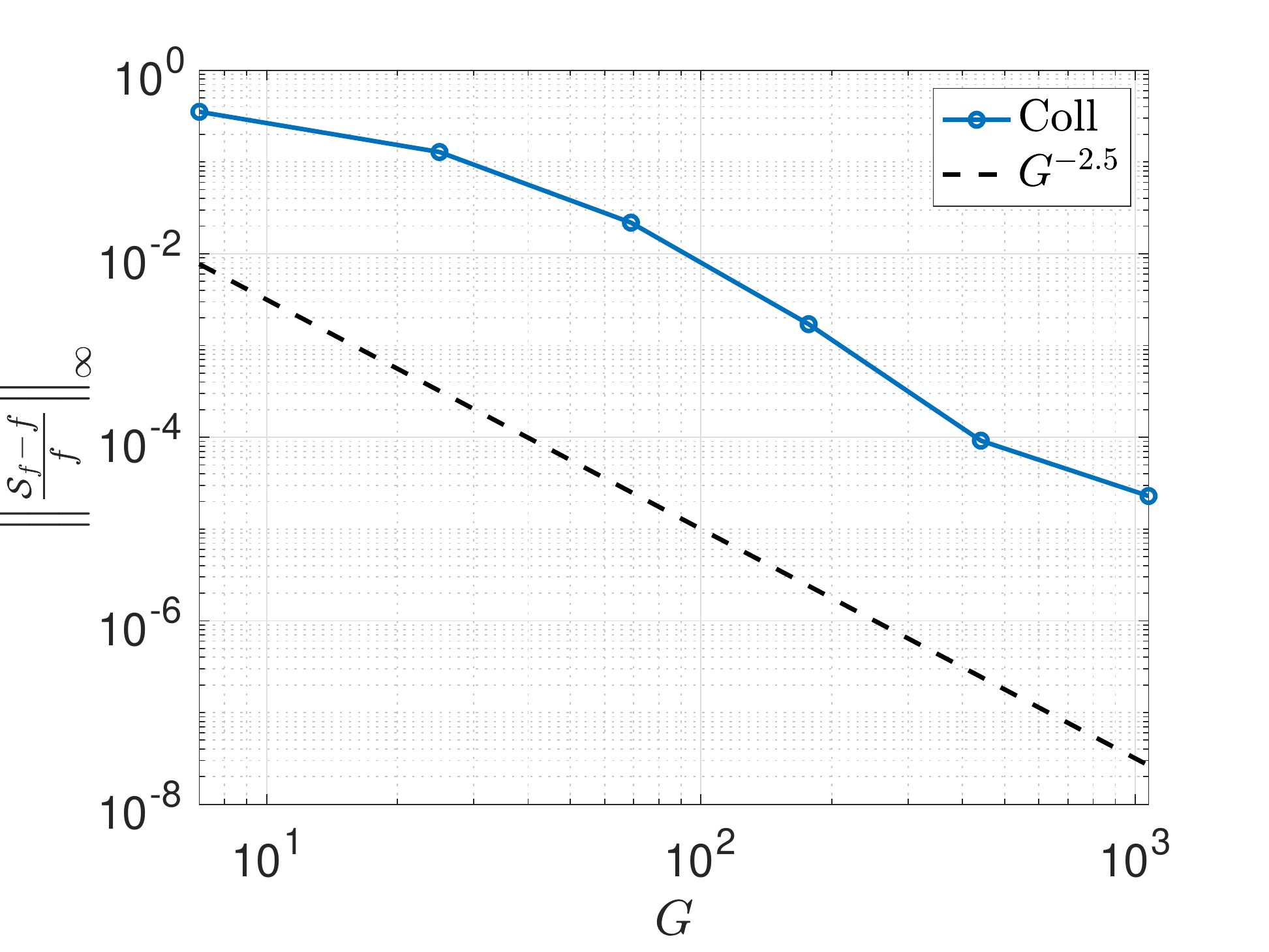}
        \label{fig:convergence_rom_1knot}}
 
     \caption{
     \ref{fig:error2_mean}  
     Relative error on the expectation of the QoI $\qoi=u(1,0,\cdot)$ plotted versus increasing cardinality $G$ of sparse grids. For the sake of comparison, the Monte Carlo error is also depicted;
     \ref{fig:convergence_rom_1knot} Maximum norm of the relative error on the QoI $\qoi=u(1,0,\cdot)$ plotted versus increasing cardinality $G$ of sparse grids.}
\end{figure*}

\begin{figure*}
    \centering
     \subfigure[]{\includegraphics[width=0.48\textwidth]{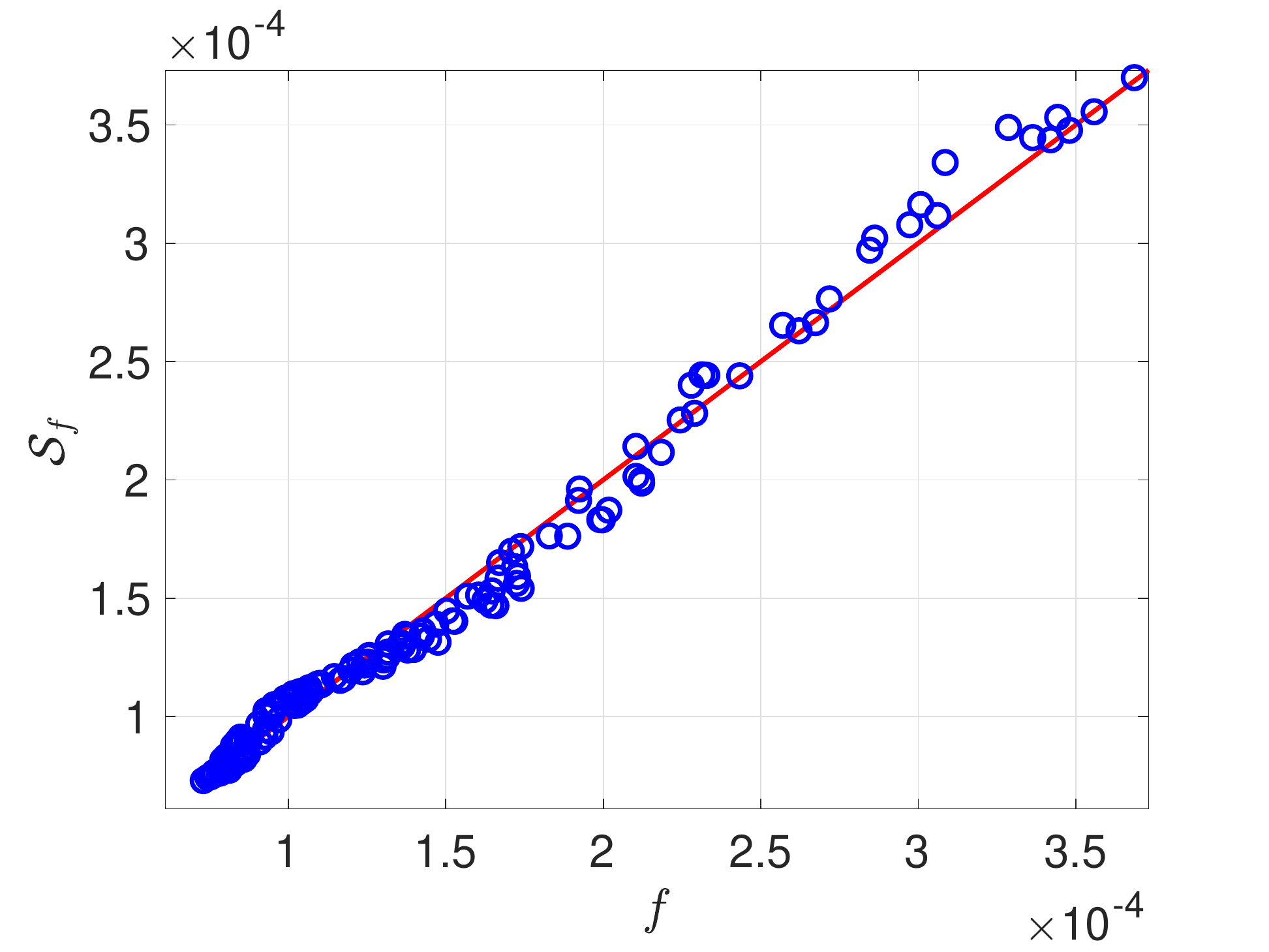}
        \label{fig:scatter_w2}}
    \subfigure[]{\includegraphics[width=0.48\textwidth]{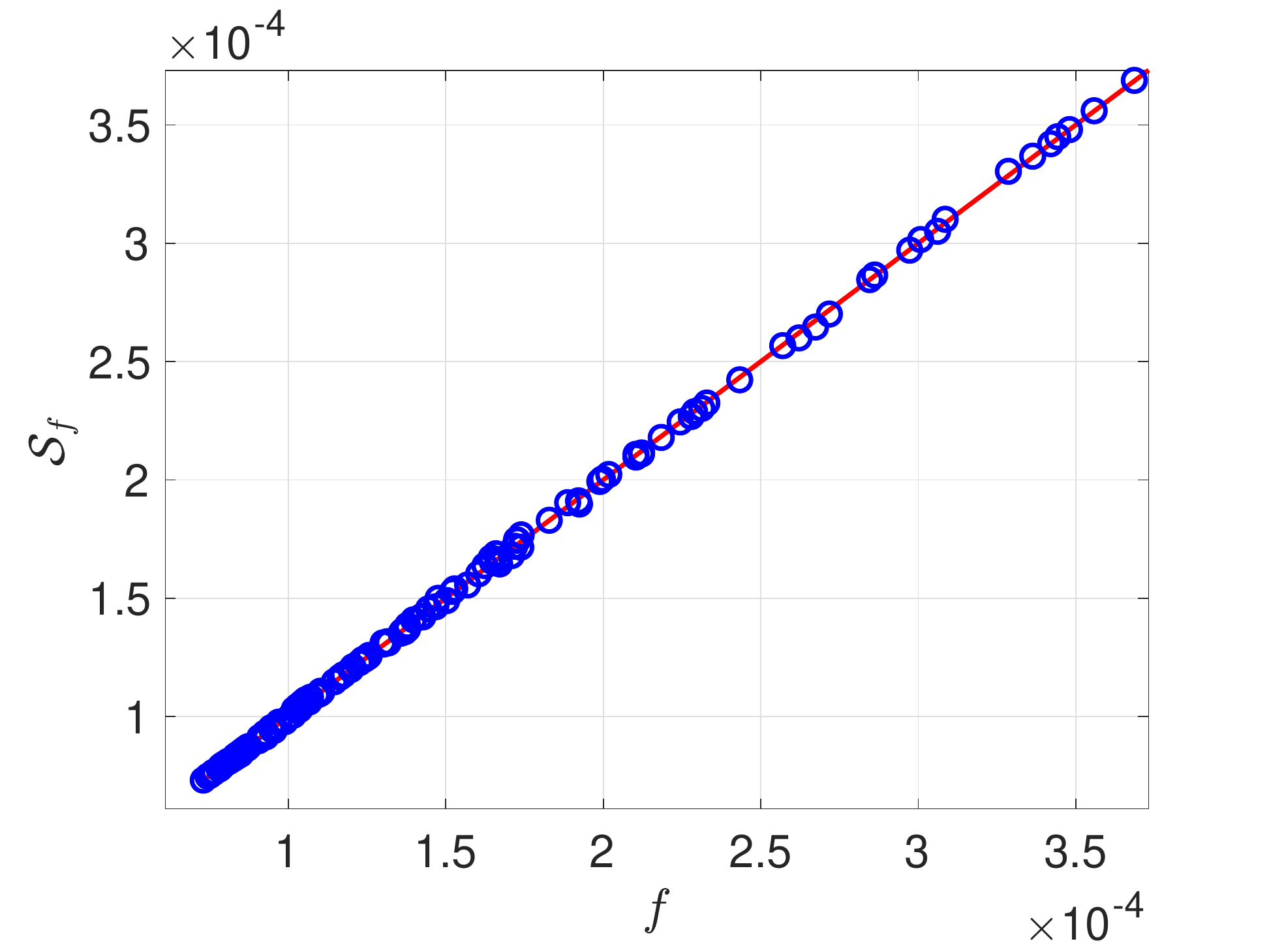}
        \label{fig:scatter_w3}}
 
     \caption{Scatterplot of the reference QoI ($x$-axis) and the surrogate QoI ($y$-axis) of level $w=2$ (Figure~\ref{fig:scatter_w2}) and level $w=3$ (Figure~\ref{fig:scatter_w3}) evaluated at the first 150 sample points $\pp^(i)$. The bisector line is depicted in red.}
\end{figure*}

\begin{figure}
    \centering
    \includegraphics[width=0.48\textwidth]{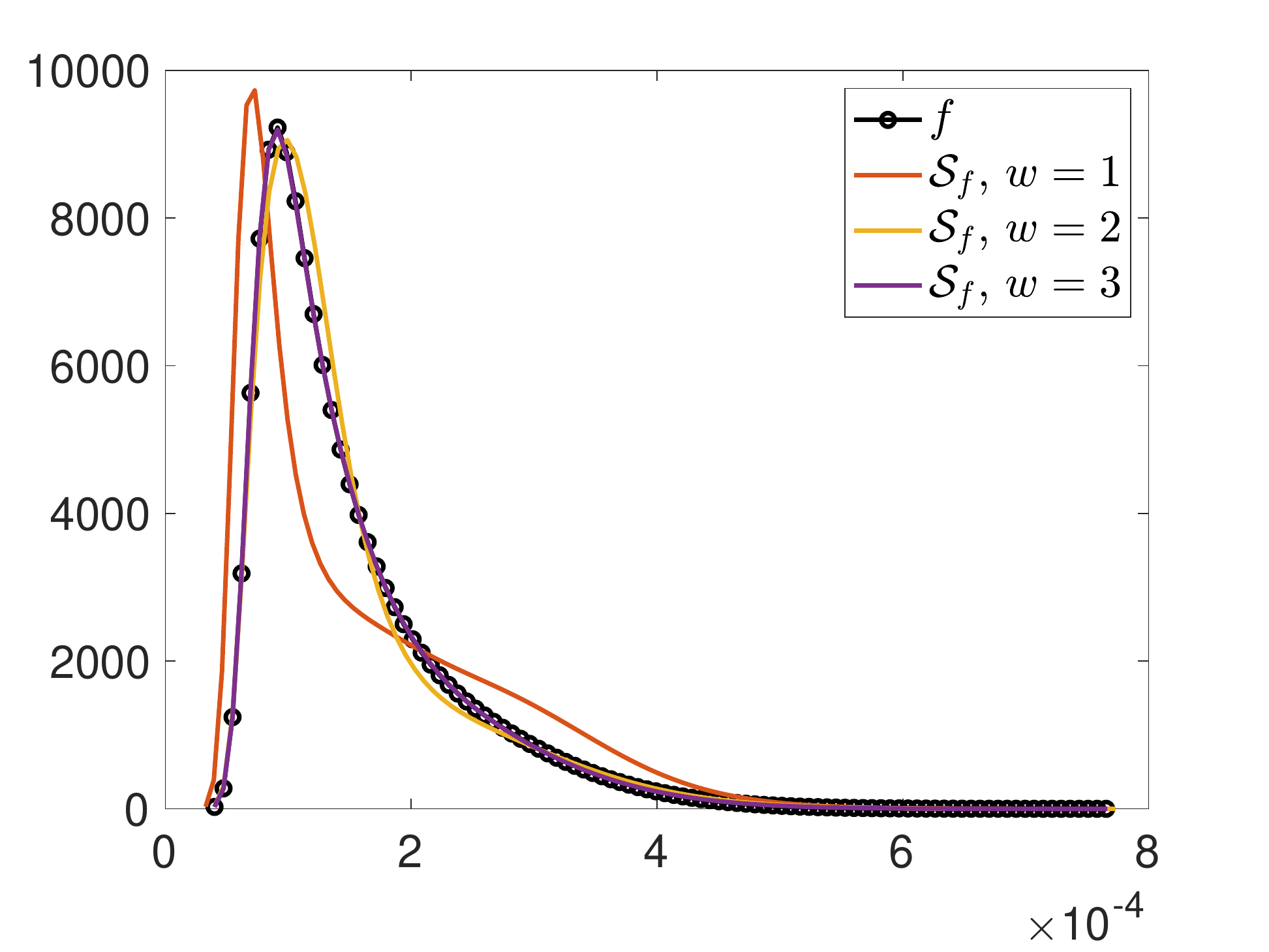}
    \caption{Approximations to the pdf of the considered $\qoi$ for increasing levels of sparse grids. The reference pdf is computed starting from evaluations of the FOM.}
    \label{fig:pdf_1knot}
\end{figure}

We now want to investigate the convergence of the sparse grid surrogate model, not only in the computation of the expected value just like we did for the previous QoI, but also in point-wise prediction. Hence, we generate $M=2000$ new samples $\{\pp^{(i)}=(p_1^{(i)},p_2^{(i)},p_3^{(i)}),\, i=1,\ldots,M)\}$ of $\pp$. For each of the new sample value, we compute the FOM solution and we compare it with the evaluations of the Smolyak sparse grid surrogate $\mathcal S_\qoi$ (see Listing~\ref{ls:evaluate}).
The relative error in the maximum norm is the given by
\begin{equation}
\label{eq:err_inf}
\left\|\frac{\mathcal S_\qoi - \qoi}{\qoi}\right\|_{\infty}
=\max_{i=1,\ldots, M}
\left|\frac{\mathcal S_\qoi(\pp^{(i)}) - \qoi(\pp^{(i)})}{\qoi(\pp^{(i)})}\right|
\end{equation}
and is displayed in Figure~\ref{fig:convergence_rom_1knot}. 
For both the expected value and the point-wise prediction, an algebraic decay of the error is observed (with estimated rate -2.7 and -2.5, respectively).
Figures~\ref{fig:scatter_w2} and~\ref{fig:scatter_w3} depict the scatterplot of the reference QoI ($x$-axis) and its surrogate ($y$-axis) of level $w=2$ and $w=3$, respectively, evaluated at the first $150$ sample points $\pp^{(i)}$ out of the 2000 samples just computed. As the level of the sparse grid increases, the blue dots tend to align along the bisector (red) line, reflecting better approximation properties of the surrogate.

Figure~\ref{fig:pdf_1knot} then graphically verifies the convergence of the pdf obtained by sampling the sparse grid surrogates $\mathcal{S}_\qoi$, for increasing levels $w=1,2,3$ (see Listing~\ref{ls:pdf}, where the built-in Matlab code \texttt{ksdensity} is used). For level $w=3$ we observe very good agreement between the reference curve and the surrogate one. Note that the alternative would be to compute the IGA solution collocated at all the $M=2000$ samples $\{\pp^{(i)}\}$, entailing a considerably larger computational effort.

\begin{lstlisting}[caption = {Matlab code to evaluate the surrogate at new parameter values.}, label={ls:evaluate}]
% generate new samples of parameter values
p1 = rand(M,1)+0.5;
p2 = rand(M,1)*0.5 + 0.25;
p3 = rand(M,1)*0.1 + 0.1;
p = [p1, p2, p3];
% point_on_grid = evaluations of QoI on the points of Sr
point_surr = interpolate_on_sparse_grid(S,Sr,point_on_grid,p');
\end{lstlisting}

\begin{lstlisting}[caption = {Matlab code to compute the pdf obtained by sampling the sparse grid surrogate of the \qoi.}, label={ls:pdf}]
% surrogate pdf for w=3 (use analogous code for w=1,2)
pdf = ksdensity(point_surr,'Support','positive');
\end{lstlisting}

\begin{lstlisting}[caption = {Matlab code to compute principal and total Sobol indices.}, label={ls:sobol}]
[Sob_princ,Sob_tot] = compute_sobol_indices_from_sparse_grid(S, Sr,f_values,domain,'legendre');
\end{lstlisting}

To conclude the UQ analysis, we compute the principal and total Sobol indices $\{S^P_i,\, i=1,2,3\}$, $\{S^T_i,\, i=1,2,3\}$ (see Section~\ref{sec:sparse_grids_for_UQ}). They are computed according to Listing~\ref{ls:sobol}, the result being
\begin{align*}
S^P&=[0.9818,0.0000,0.0088],\\
S^T&=[0.9912,0.0001,0.0182].
\end{align*}
They confirm that the variability of the second parameter $p_2$ does not affect the surrogate value, as was previously observed by means of Figure~\ref{fig:surface}, and moreover they hint that the third parameter
plays a negligible role as well.

As a conclusion of the analysis carried out, we can state that 
thanks to the sparse grids machinery, a small computational effort
was enough to carry out a UQ analysis that allows to draw these conclusions
on the timeber model at hand:
(i) the parameter playing the most important role in the model is $p_1$; (ii) small variability of the QoI is caused by $p_3$; (iii) $p_2$ affects the QoI in a negligible way. These results are as expected, since $p_2$ and $p_3$ have local effects on the solution to the PDE~\eqref{eq:pde}, whereas the considered QoI is affected by global quantities, only.

% % % % % % % % % % % % % % % % % % % % 
\subsection{Two-knots example}
\label{sec:two_knot}
% % % % % % % % % % % % % % % % % % % % 

\begin{figure*}
\centering
    \subfigure[First sample of $E(x,y,\pp)$]
    {\includegraphics[width=0.455\textwidth]{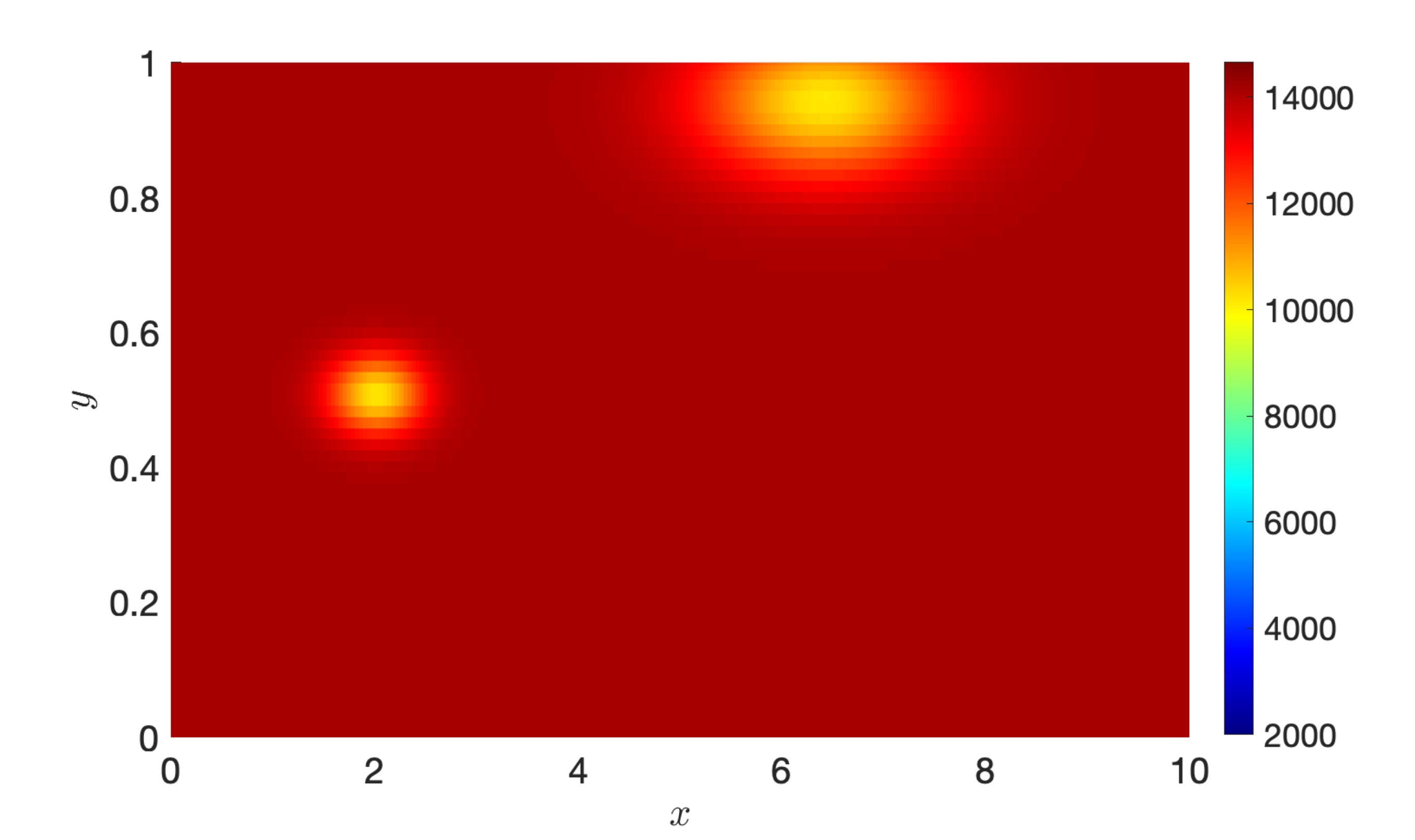}}
    \subfigure[First solution $u_x(x,y,\pp)$]
    {\includegraphics[width=0.455\textwidth]{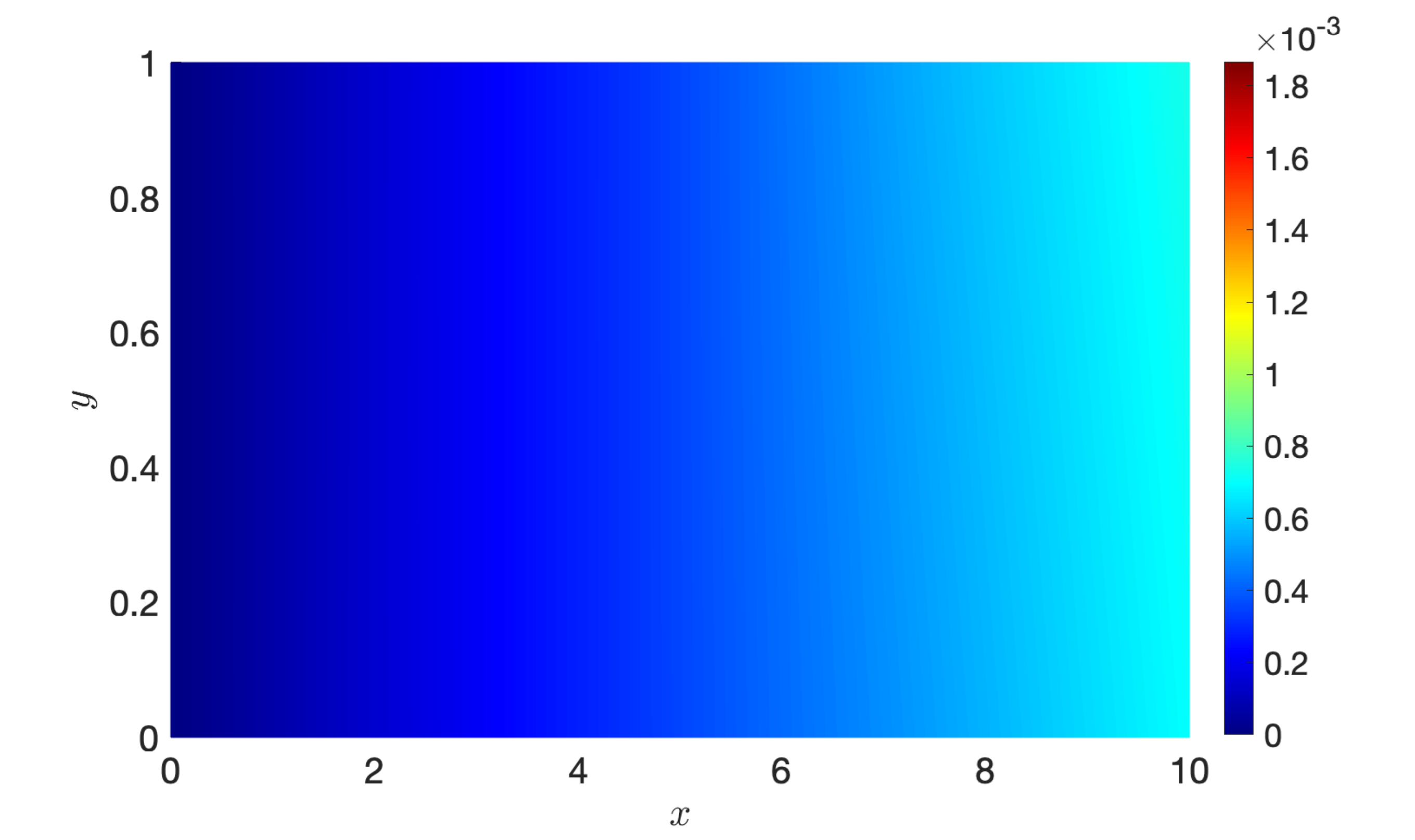}} 
    
    \subfigure[Second sample of $E(x,y,\pp)$]
    {\includegraphics[width=0.455\textwidth]{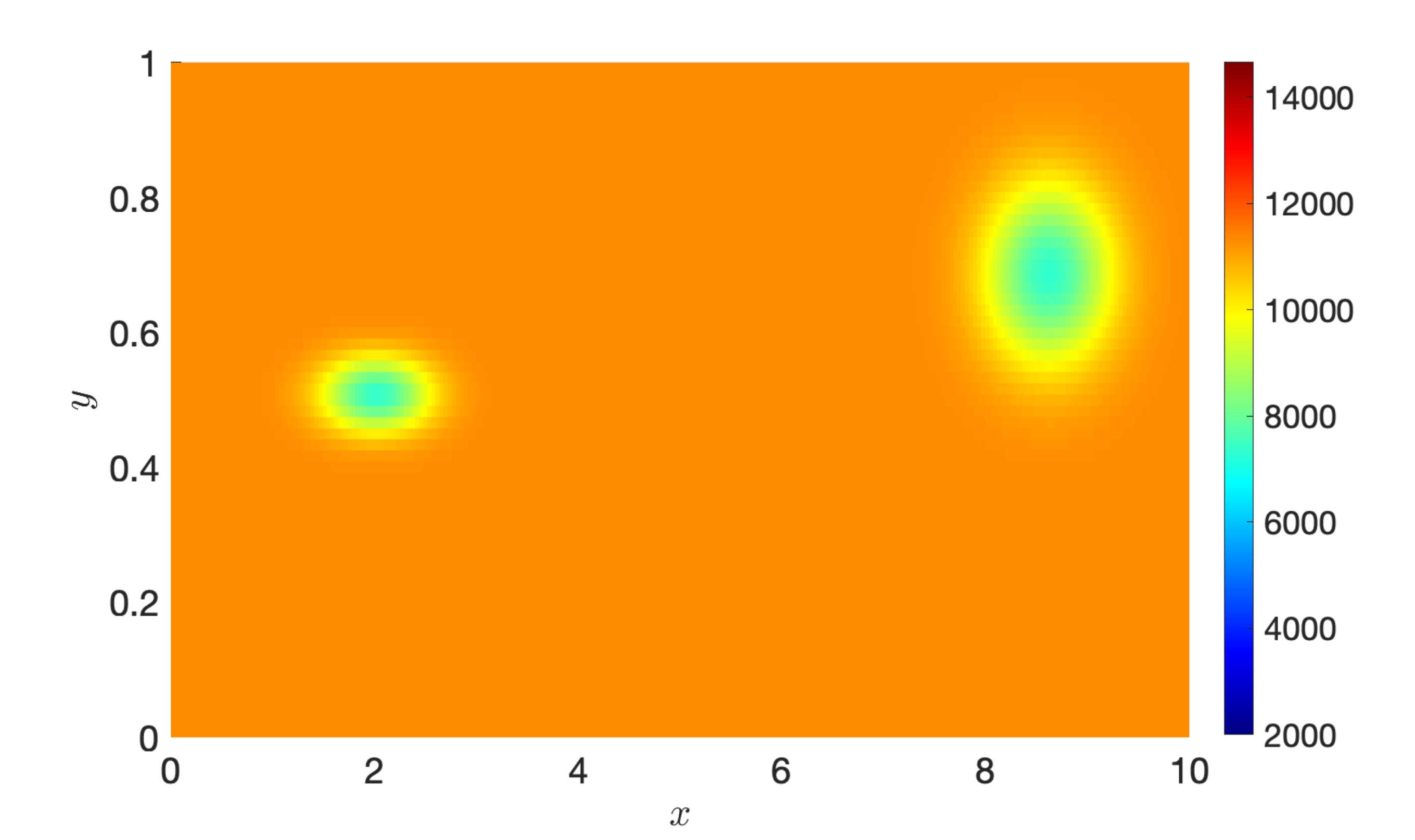}}
    \subfigure[Second solution $u_x(x,y,\pp)$]
    {\includegraphics[width=0.455\textwidth]{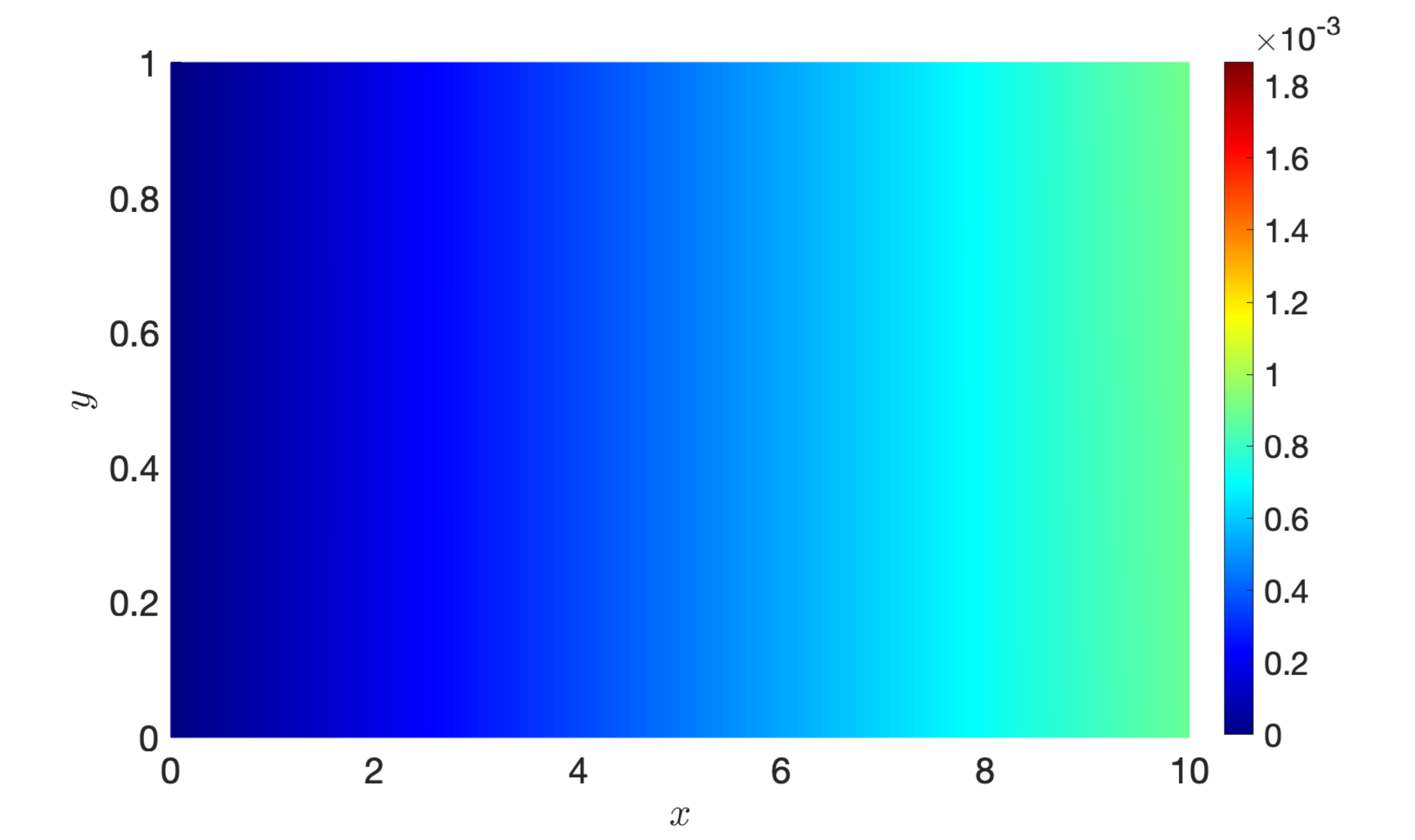}} 
    
    \subfigure[Third sample of $E(x,y,\pp)$]
    {\includegraphics[width=0.455\textwidth]{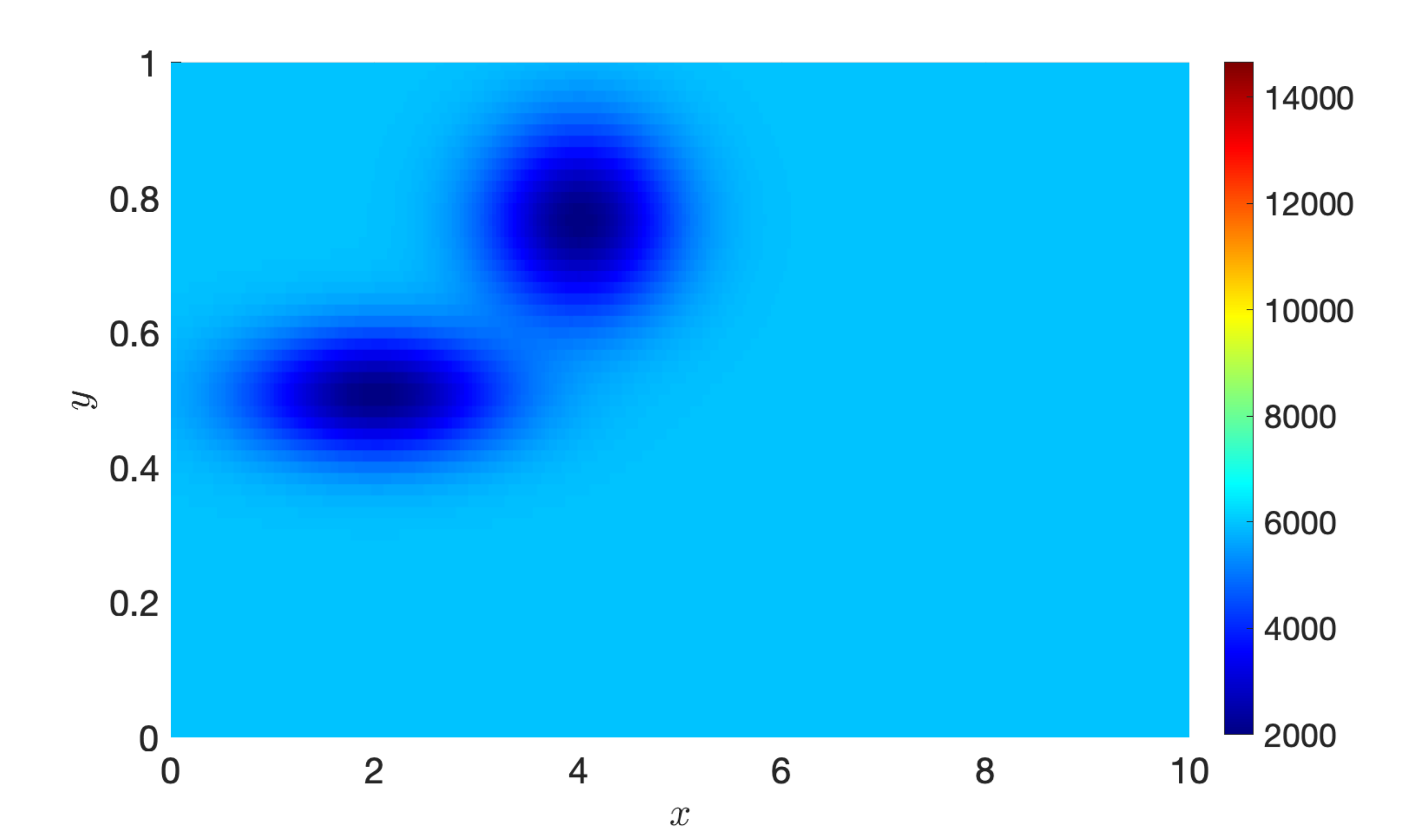}}
    \subfigure[Third solution $u_x(x,y,\pp)$]
    {\includegraphics[width=0.455\textwidth]{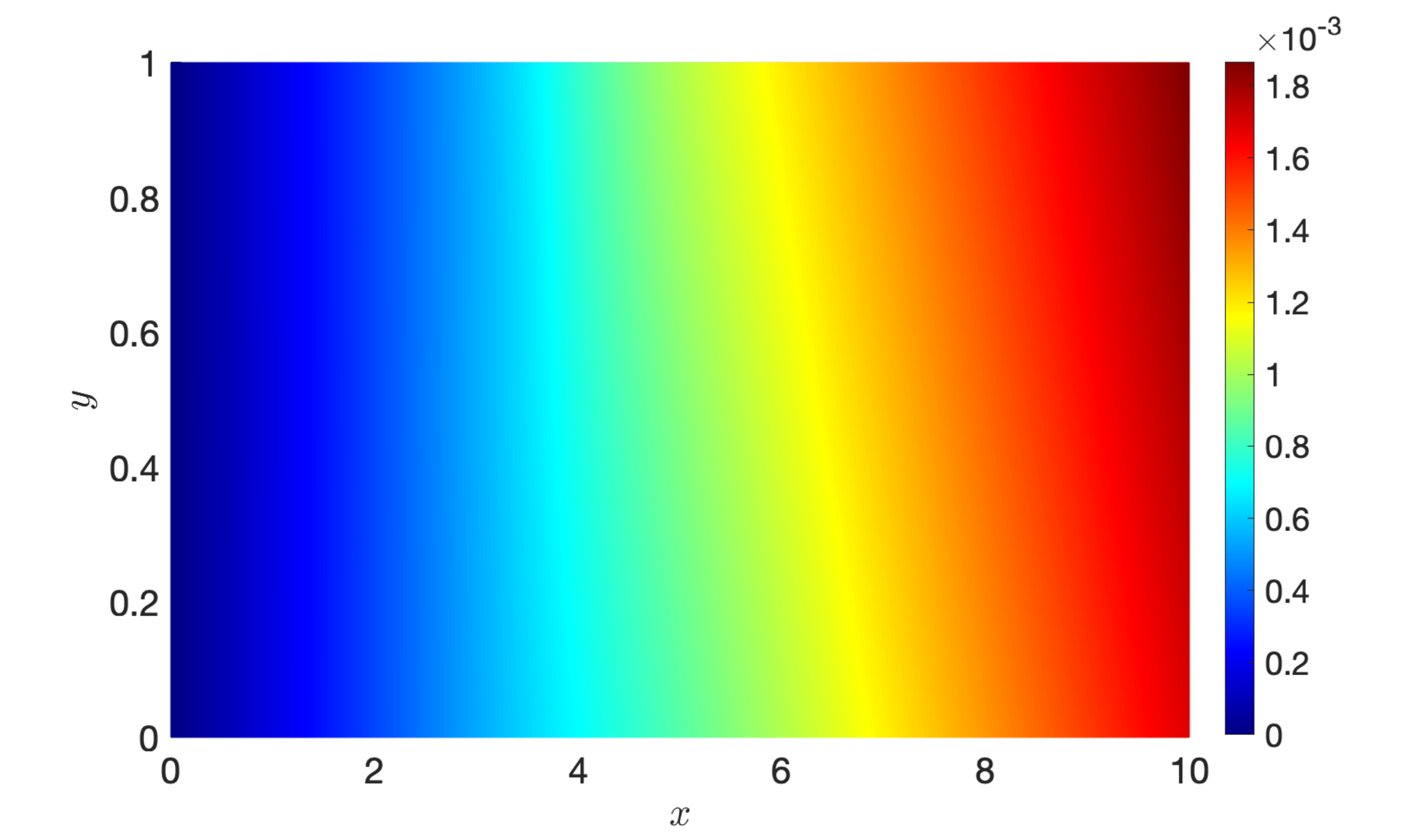}} 
    
    \subfigure[Fourth sample of $E(x,y,\pp)$]
    {\includegraphics[width=0.455\textwidth]{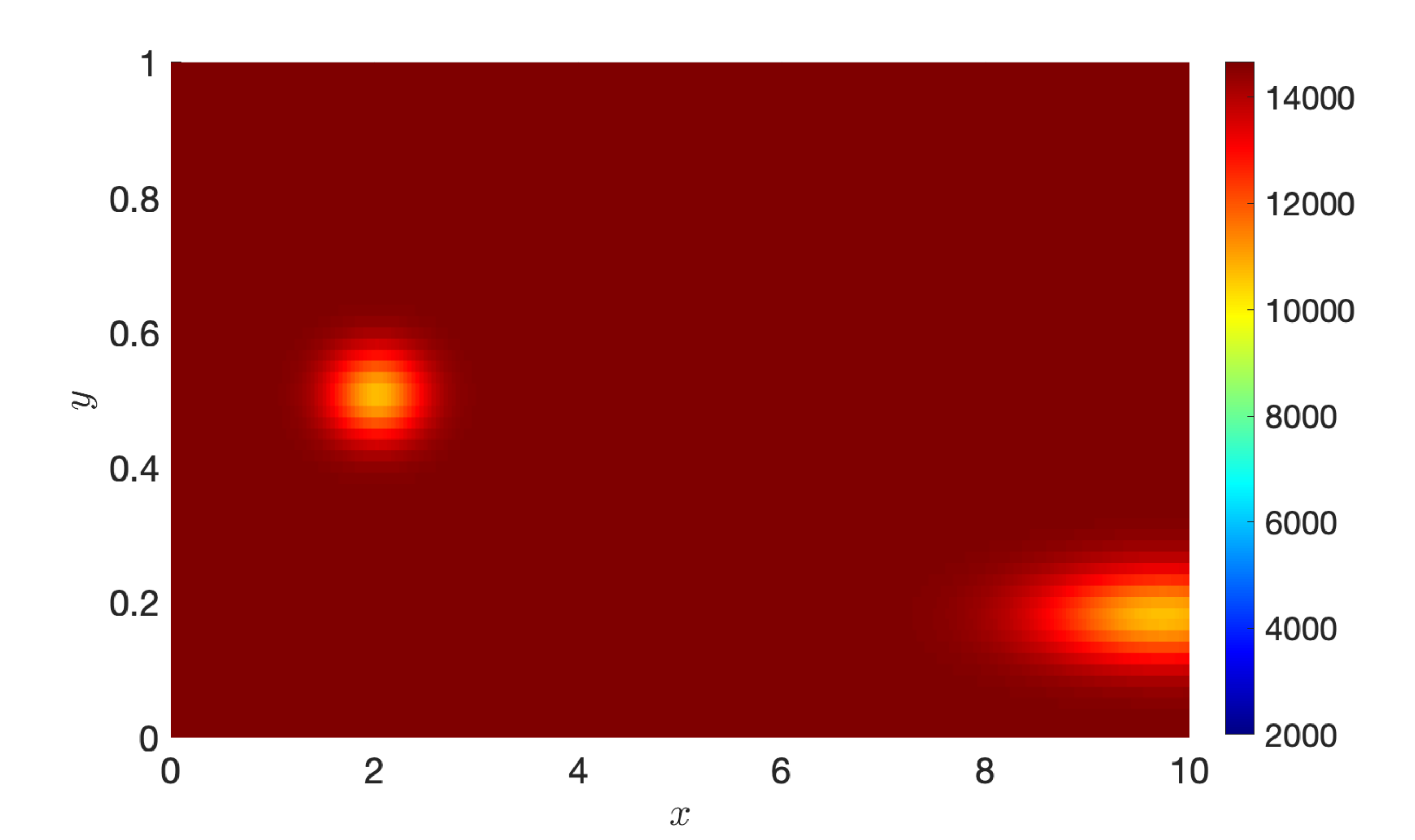}}
    \subfigure[Fourth solution $u_x(x,y,\pp)$]
    {\includegraphics[width=0.455\textwidth]{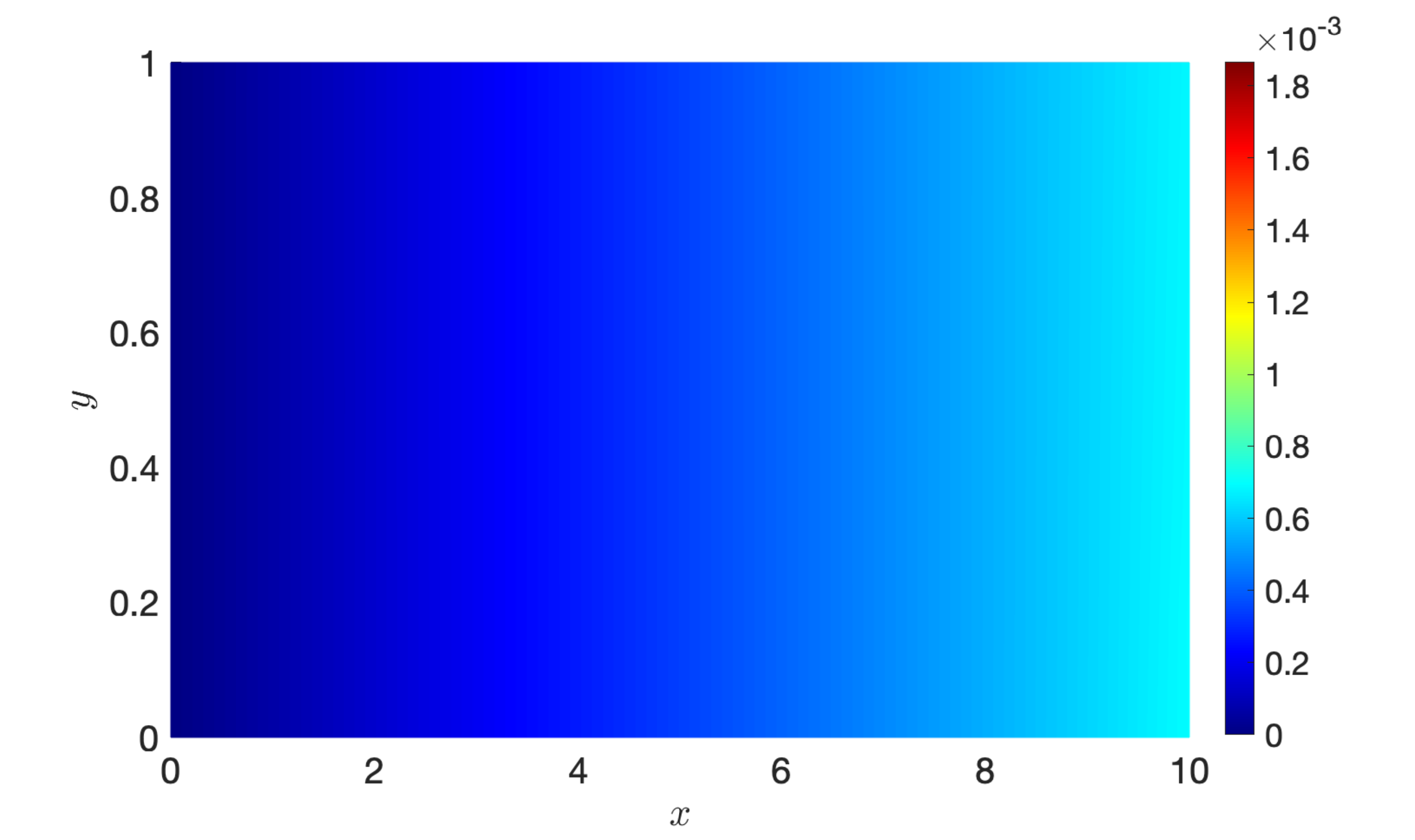}} 
    
    \caption{Plot of four samples of $E(x,y,\pp)=E_0\alpha(x,y,\pp)$, with $\alpha(x,y,\pp)$ as in~\eqref{eq:alpha_2knot} and of the corresponding horizontal displacement $u_x(x,y,\pp)$ computed via the IGA method (Section~\ref{sec:IGA}) for $\pp=(14134,4.3976, 0.4340, 0.3323, 0.8566,0.0414,0.0789)$, $\pp=(11324, 6.6020, 0.1787, 0.3680,0.4308,0.0383, 0.0924)$, $\pp=(5.9754, 1.9932, 0.2577,0.8764,0.6428,0.0649,0.0972)$ and $\pp=(1.4649,7.7164, -0.3288,0.3241,0.8283,0.0457,0.0480)$ (from top to bottom).}
    
    \label{fig:2knot_samples}
    
\end{figure*}

Let us take $L = 10 \meter \mbox{ and } H = 1 \meter$, and choose 
\begin{align}
    \nonumber
    &\alpha(x,y,\pp)
    = p_1 - \gamma_1 \exp\left(-\frac{(x-\bar x)^2}{2p_2^2}\right)
    \exp\left(-\frac{(y-\bar y)^2}{2p_4^2}\right)\\
    \label{eq:alpha_2knot}
    &\quad - \gamma_2 \exp\left(-\frac{(x-\bar x-p_6)^2}{2p_3^2}\right)
    \exp\left(-\frac{(y-\bar y-p_7)^2}{2p_5^2}\right),
\end{align}
with parameter $\pp=(p_1,p_2,p_3,p_4,p_5,p_6,p_7)$ whose entries are $p_1\sim\mathcal{U}(0.5,1.5)$, $p_2,p_3\sim\mathcal{U}(0.3,1)$, $p_4,p_5\sim\mathcal U(0.03,0.1)$, $p_6\sim\mathcal{U}(1,8)$ and $p_7\sim\mathcal{U}(-0.5,0.5)$, and with fixed values $\gamma_1=\gamma_2=0.4$.
In particular, in this second example, we aim at modeling the presence of two knots along the timber beam: one with fixed coordinates $(\bar x,\bar y)$ and the second at random distance from the first one with coordinates $(\bar x+p_6,\bar y+p_7)$; the value of $E_0$ is again $E_0 = 10^4 \mega\pascal$.
As in the first example, $p_1$ - when multiplied by $E_0$ - controls the nominal value of the Young modulus away from the knot. The remaining parameters model the width of the two knots along the horizontal ($p_2,\, p_3$) and vertical ($p_4,\, p_5$) direction.
Figure~\ref{fig:2knot_samples} depicts four samples of $E(x,y,\pp)=E_0\alpha(x,y,\pp)$, with $\alpha(x,y,\pp)$ as in~\eqref{eq:alpha_2knot} (left column) and the first component of the corresponding IGA solution (right column). 
We underline that in the considered setting the second knot can be placed (i) close to the top/bottom boundary of the beam (first sample, showing the case of proximity to the top boundary); (ii) well-contained in the beam and distant from the first knot (second sample); (iii) well-contained in the beam but close to the first knot (third sample); (iv) close to the right boundary of the beam (fourth sample). 

In this second example, we are dealing with a problem with a larger number of uncertain parameters, therefore we consider not only the classical Smolyak sparse grids used in the previous example, but also the more effective a-posteriori adaptive sparse grids. In this version of sparse grids, multi-indices $\ii$ are added to the multi-index set $\II$ in Equation \eqref{eq:sparsegrid-combitec} in an iterative way, following a simple yet powerful procedure based on an error-cost criterion (see e.g. \cite{gerstner.griebel:adaptive,nobile.eal:adaptive-lognormal} for details): 
\begin{itemize}
    \item a number of potential candidates $\jj$ is added to the sparse grid
     to the multi-index set $\II$;
    \item for each of them a profit indicator 
    is computed (i.e. the ratio between the change in the prediction of
    $\Ev[\qoi]$ due to having added $\jj$ to the sparse grid 
    and the number of new FOM evaluations requested by it);
    \item the candidate with the largest profit is selected and added to 
    $\II$, and the set of candidates is updated accordingly.
\end{itemize}
This algorithm is usually very effective in quickly
determining a good set $\II$, although it is not entirely optimal
in terms of cost since the profits are evaluated 
only \emph{after} having performed the corresponding FOM evaluations 
(hence the name ``a-posteriori adaptive''), therefore some
computational cost is ``wasted'' to detect multi-indices with
small profit.

\begin{figure}
    \centering
    \includegraphics[height=0.3\textwidth]{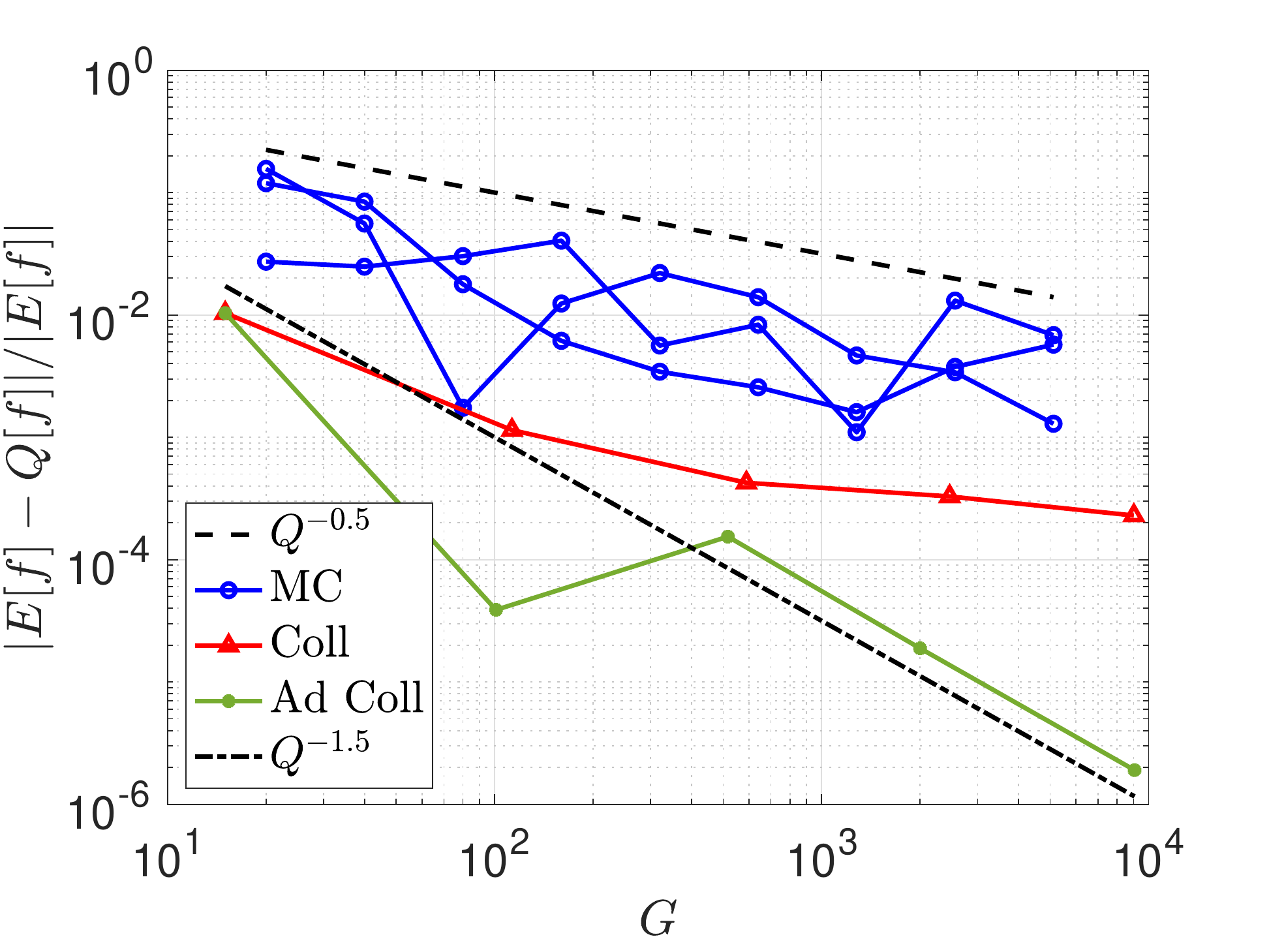}
     \caption{
     Error convergence of the sparse grid surrogates.
     As a comparison, three instances of Monte Carlo errors are also reported. The lines are plotted versus the number of PDE solves (the cardinality $G$ of the sparse grid the Collocation method, the number of samples for the Monte Carlo method).}
    \label{fig:error_point}
\end{figure}

Coming back to the computational example,
let us consider the same real-valued QoI as in the first example, namely the evaluation of the horizontal displacement at the bottom-right corner of the beam $\qoi(\pp)=u_x(10,0,\pp)$.
The reference value $\Ev[u_x(10,0,\cdot)]$ is approximated  using an a-posteriori adaptive sparse grid with 30105 collocation points (see Listing~\ref{ls:adaptive}) and it is compared with its approximation $\mathcal Q[u_x(10,0,\cdot)]$ computed either using Smolyak sparse grids of increasing level $w=1,\ldots,6$ (red line in Figure~\ref{fig:error_point}) or a-posteriori adaptive sparse grids with increasing number of collocation points (magenta line in  Figure~\ref{fig:error_point}). 

\begin{figure}
    \centering
    \includegraphics[height=0.3\textwidth]{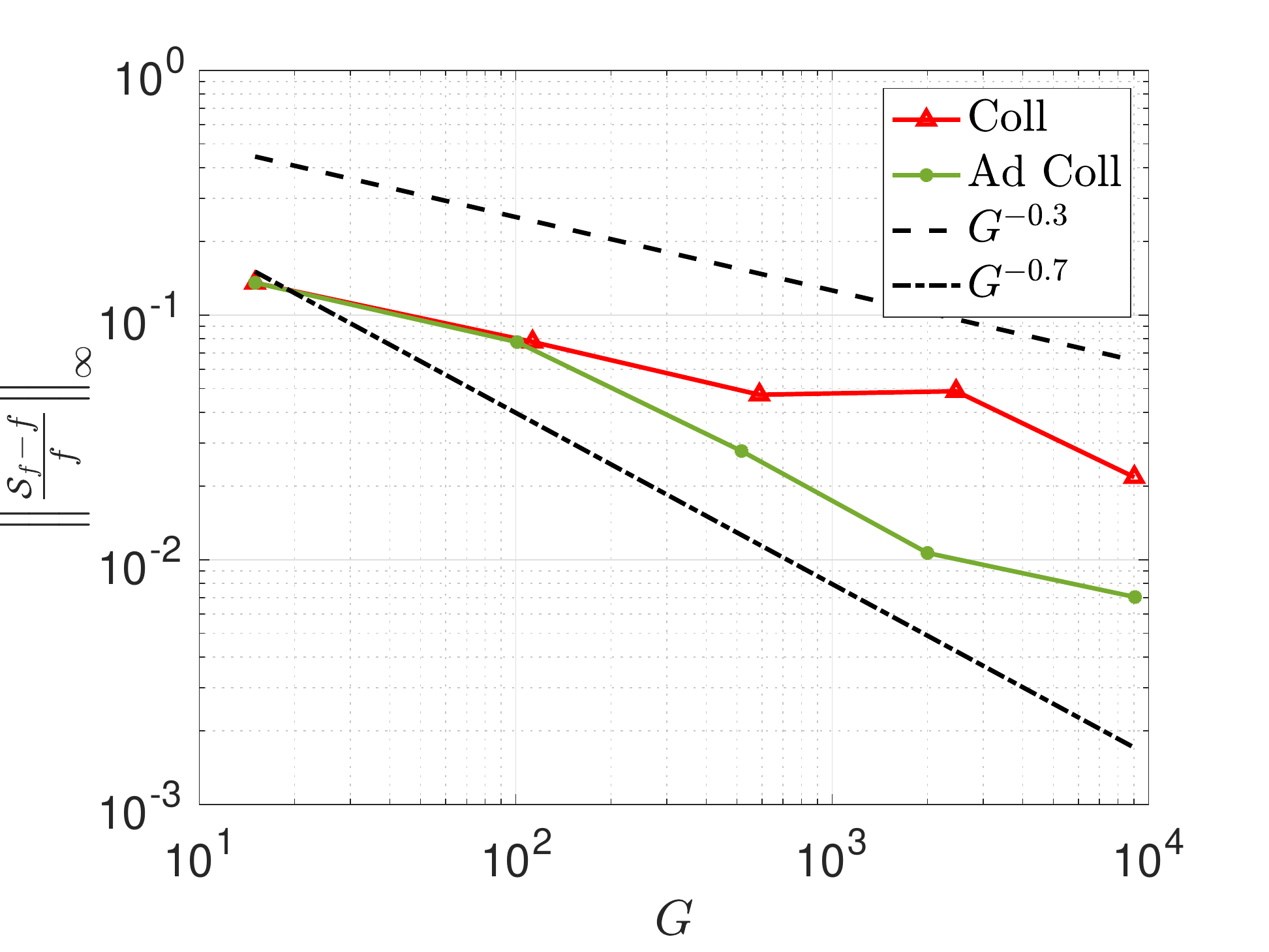}
     \caption{
     Maximum norm of the relative error on the QoI $\qoi=u_x(10,0,\cdot)$ plotted versus increasing cardinality of sparse grids.}
    \label{fig:convergence_rom_2knots}
\end{figure}

Following the same lead as in the first example, we now want to investigate the convergence of the sparse grids surrogate to the FOM. 
Hence, the QoI is computed by the FOM at $M=5000$ randomly generated samples of $\pp$ (denoted as $\{\pp^{(i)},\, i=1,\ldots,M\}$. 
In contrast, the a-posteriori adaptive and the Smolyak sparse grid surrogates are evaluated at all points $\pp^{(i)}$ and the largest relative error is computed by Equation \eqref{eq:err_inf}. 
The decay of this approximation error
as the sparse grids construction cost increases 
is depicted in Figure~\ref{fig:convergence_rom_2knots}. 
Both Figure~\ref{fig:error_mean} and~\ref{fig:convergence_rom_2knots} display improved rates of convergence of the a-posteriori adaptive sparse grids, when compared to the non-adaptive ones (i.e., the Smolyak sparse grids).

\begin{figure}
    \centering
    \includegraphics[height=0.3\textwidth]{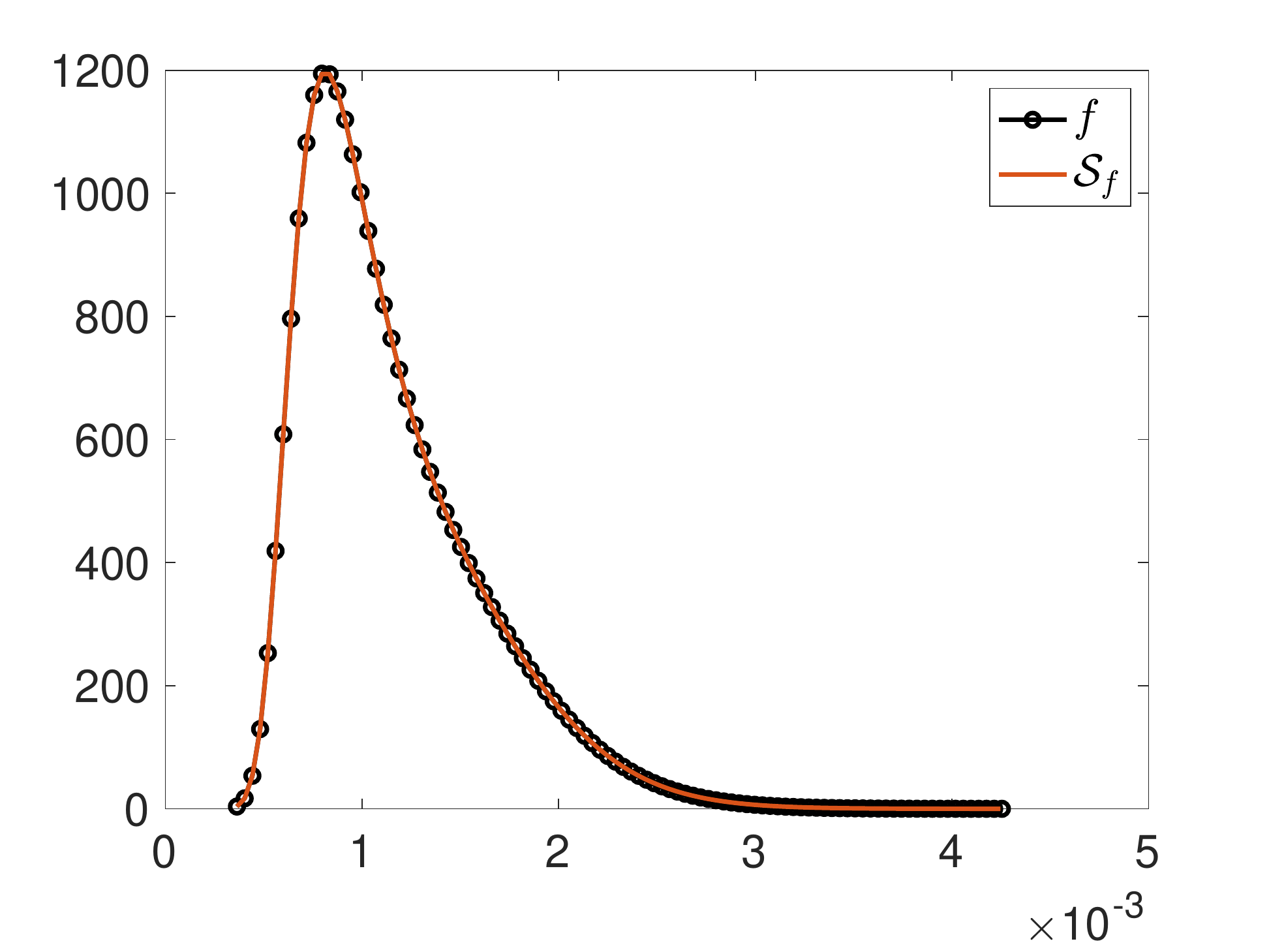}
     \caption{
     Approximation to the pdf of the considered $\qoi$ corresponding to the surrogate obtained by means of the a-posteriori adaptive sparse grid. The reference pdf is computed starting from evaluations of the FOM.}
    \label{fig:pdf_2knots}
\end{figure}

Next, we graphically verify the convergence of the pdf obtained by sampling the a-posteriori adaptive sparse grid surrogate,
see Figure \ref{fig:pdf_2knots}: a very good agreement between
he exact and the surrogate pdfs can be observed. 
Finally, the principal and total Sobol indices are computed: 
\begin{align*}
S^P&=[0.9914,0.0000,0.0040,0.0001,0.0001,0.0002,0.0002],\\
S^T&=[0.9945,0.0002,0.0075,0.0002,0.0008,0.0003,0.0010].
\end{align*}
As observed in the first example, the first parameter $p_1$ is by far the most important one, in the sense that it essentially affects all the variability of the selected QoI. The other parameters play a much smaller role, and in particular the second one $p_2$ (horizontal width of the first knot) appears to be negligible.

\begin{lstlisting}[caption = {Matlab code to compute an adaptive sparse grid approximation of the QoI.}, label={ls:adaptive}]
% number of parameters
N=7; 
% knots for p1, p2, p3, p4, p5, p6 and p7
knots_Y0=@(n) knots_CC(n,0.5,1.5,'prob'); % c
knots_Y1=@(n) knots_CC(n,1,8,'prob'); % r
knots_Y2=@(n) knots_CC(n,-0.5,0.5,'prob'); % s
knots_Y3=@(n) knots_CC(n,0.3,1,'prob'); % sx1
knots_Y4=@(n) knots_CC(n,0.3,1,'prob'); % sx2
knots_Y5=@(n) knots_CC(n,0.03,0.1,'prob'); % sy1
knots_Y6=@(n) knots_CC(n,0.03,0.1,'prob'); % sy2
knots = {knots_Y0,knots_Y1,knots_Y2,knots_Y3,knots_Y4,knots_Y5,knots_Y6};
[lev2knots,idxset]=define_functions_for_rule('SM',N);
% number of maximum collocation points 
Max_Points = 30000;
% QoI (to be implemented separately)
f = @(y) QoI(y);
adapt = adapt_sparse_grid(f,N,knots,lev2knots)
\end{lstlisting}

% % % % % % % % % % % % % % % % % % % % 
\section{Conclusions}
\label{sec:conclusions}
% % % % % % % % % % % % % % % % % % % % 

In this paper we used the Sparse Grid Matlab kit (freely available on-line) for the UQ of the displacements in the field of continuum linear mechanics.
The considered problem has been discretized by means of the IGA collocation (in the space variables), while the dependence on random parameters has been treated by means of the stochastic collocation method on both Smolyak and a-posteriori adaptive sparse grids.

The use of the Sparse Grid Matlab kit presents several advantages. First, it can be easily interfaced with any black-box solvers for the (deterministic) mechanical problem. Second, the numerical methods implemented in the kit outperform standard UQ techniques, like the plain Monte Carlo method. Finally, it provides outputs that can be readily interpreted and exploited in the engineering practice.

Future work will include the analysis of more sophisticated and possibly anysotropic mechanical problems, where the variability of grain direction is accounted for.

% % % % % % % % % % % % % % % % % % % % 
\section{Declaration of Competing Interest}
% % % % % % % % % % % % % % % % % % % % 

The authors declare that they have no known competing financial interests or personal relationships that could have appeared to influence the work reported in this paper.

% % % % % % % % % % % % % % % % % % % % 
\section{Acknowledgments}
% % % % % % % % % % % % % % % % % % % % 

Lorenzo Tamellini hs been supported by the PRIN 2017 project 201752HKH8 ``Numerical Analysis for Full and Reduced Order Methods for the efficient and accurate solution of complex systems governed by Partial Differential Equations (NA-FROM-PDEs)''.
Lorenzo Tamellini also acknowledges the support of GNCS-INdAM (Gruppo Nazionale Calcolo Scientifico - Istituto Nazionale di Alta Matematica), Italy.
F. Bonizzoni is member of the INdAM Research group GNCS and her work was part of a project that has received funding from the European Research Council
ERC under the European Union's Horizon 2020 research and innovation program (Grant
agreement No. 865751).
Finally, the authors would like to acknowledge Alessandro Reali for his contribution in the development of the IGA collocation code.

\appendix

% % % % % % % % % % % % % % % % % % % % 
\section{Basics on B-splines}
\label{appendix:basic_b-splines}
% % % % % % % % % % % % % % % % % % % % 

Let us introduce two knot vectors:
\begin{equation}
\begin{aligned}
& X=\{x_1=0\leq x_1\leq \dots\leq x_{\ncp+\degr+1}=L\}
\\
& Y=\{y_1=0\leq y_1\leq \dots\leq y_{\mcp+\degq+1}=H\}  
\end{aligned}
\end{equation}

where $\degr$ and $\degq$ are the degree of the B-splines and $\ncp$ and $\mcp$ are the numbers of basis functions. Pairs $(x_i,y_j)\in X\times Y$ correspond to coordinates of points in the 2D domain $D$. In particular, we take $X, Y$ as so-called \emph{open} vectors, i.e., the first and last knots of $X$ ($Y$, respectively) have multiplicity $\degr+1$ ($\degq+1$, respectively), and - for simplicity - we choose uniformly equispaced $X$ and $Y$ knots.

Given the knot vector $X$, the uni-variate B-spline basis functions in the $x$-variable are defined recursively as follows:
\begin{itemize}
    \item for $\degr=0$:
    \begin{equation*}
        N_{i}^{0}(x)=
        \begin{cases}
            1, & \mbox{ if } x_i\leq x < x_{i+1},\\
            0, & \mbox{ otherwise},
        \end{cases}
    \end{equation*}
    \item for $\degr>1$:
    \begin{equation*}
        N_{i}^{\degr}(x)=
        \begin{cases}
            \cfrac{x-x_i}{x_{i+\degr}-x_i}N_{i,\degr-1}(x)+ 
            \cfrac{x_{i+\degr+1}-x} {x_{i+\degr+1}-x_{i+1}}N_{i+1,\degr-1}(x),  \qquad
            \\     
            \hfill \mbox{if }  x_i\leq x < x_{i+\degr+1},    
            \\
            0, \hfill\mbox{ otherwise},
        \end{cases}
    \end{equation*}
\end{itemize}
with the convention $0/0 = 0$.
Given the knot vector $Y$, the uni-variate B-spline basis functions $M_{j}^{\degq} \left( y \right)$ are defined analogously. 
The tensor product construction leads to bi-variate basis functions for the 2D domain $D$, given by
\begin{equation*}
R_{i,j}^{\degr,\degq} \left( x,y \right) = N_{i}^{\degr} \left( x \right) M_{j}^{\degq} \left( y \right)
\end{equation*} 
for all $i=1,\ldots,\ncp$, $j=1,\ldots,\mcp$.

% % % % % % % % % % % % % % % % % % % % 
\section{Formulas for sparse grids} \label{appendix:sparse_grids_formulas}
% % % % % % % % % % % % % % % % % % % % 

In this appendix, we report some auxiliary formulas for readers who are interested 
in the details of how the tensor interpolant $\qoi_{m(\ii)}$ in Equation \eqref{eq:sparsegrid-combitec} is obtained. To this end, let us recall that $\ii$ is a vector of $N$ integers, and that $m(\cdot)$
is an increasing function. We then let:
\begin{itemize}
    \item $\GGG_{n,m(i_n)}$ a set of $m(i_n)$ points in the range of the $n$-th parameter $p_n$,
    (for instance, the Clenshaw--Curtis points introduced in Equation~\eqref{eq:points}), 
    \[
    \GGG_{n,m(i_n)} = \{p_{n,m(i_n)}^1,p_{n,m(i_n)}^2,\ldots,p_{n,m(i_n)}^{m(i_n)}\}
    \]
    \item $\ell_{n,m(i_n)}^{k}$ be the Lagrange polynomial associated to the $k$-th node of $\GGG_{n,m(i_n)}$, i.e., a polynomial that has value 1 in $p_{n,m(i_n)}^{k}$ and 0 in 
    every other point of $\GGG_{n,m(i_n)}$. The explicit expression of $\ell_{n,m(i_n)}^{k}(p)$ reads:
    \[
    \ell_{n,m(i_n)}^{k}(p) = 
    \prod_{j=1, j \neq k}^{m(i_n)} \frac{p - p_{n,m(i_n)}^{j}}{p_{n,m(i_n)}^{k}-p_{n,m(i_n)}^{j}}.
    \]
    \item $\GGG_{m(\ii)}$ is the cartesian product of the univariate sets $\GGG_{n,m(i_n)}$, for $n=1,\ldots,N$, namely
    \[
    \GGG_{m(\ii)} = \prod_{n=1}^N \GGG_{n,m(i_n)}.
    \]
    It contains $m(i_1) \times m(i_2) \times \cdot m(i_N)$ points. Each of these points corresponds a multi-index $\jj$, that is component-wise smaller than $m(\ii)$, i.e.,
    \[
    \GGG_{m(\ii)} = \{ \pp_{\jj} \in \Rset^N : p_{j_n} = x_{n,m(i_n)}^{j_n} \mbox{with } j_n \leq m(i_n) \}.    
    \]
    \item To each $\pp_{\jj}$ we can associate the multi-variate Lagrange polynomial given by
    \[
    \ell_{m(\ii)}^{\jj}(\pp) = \prod_{n=1}^N \ell_{n,m(i_n)}^{j_n}(p_n).    
    \]
\end{itemize}
With these definitions in place, we can finally define the tensor interpolant $\qoi_{m(\ii)}$ as
\[
\qoi_{m(\ii)}(\pp) = \sum_{\pp_{\jj} \in \GGG_{m(\ii)}} \qoi(\pp_{\jj}) \ell_{m(\ii)}^{\jj}(\pp).
\]
The sparse grid $\GGG$ (i.e., the 
union of the points required to assemble each $\qoi_{m(\ii)}$ in Equation~\eqref{eq:sparsegrid-combitec}) can
be obtained as
\[
\GGG = \bigcup_{\ii \in \II} \GGG_{m(\ii)}.
\]

\end{document}